\documentclass[10pt,twoside,reqno]{amsart}
\usepackage{latexsym, amsmath, amssymb, longtable, booktabs,amscd}
\usepackage{amsmath,amssymb,amsthm}
\usepackage[overload]{empheq}
\usepackage{xcolor}
\usepackage{enumitem}
\usepackage[all]{xy}
\UseComputerModernTips
\numberwithin{equation}{section}
\usepackage[hidelinks]{hyperref} 
\newtheorem{thm}{Theorem}[section]
\newtheorem{theorem}[thm]{Theorem}
\newtheorem{cor}[thm]{Corollary}
\newtheorem{conj}[thm]{Conjecture}
\newtheorem{prop}[thm]{Proposition}
\newtheorem{lemma}[thm]{Lemma}
\theoremstyle{definition}

\theoremstyle{remark}
\newtheorem{rk}[thm]{Remark}

\newcommand{\kerp}{\mathrm{Ker} (P)}
\newcommand{\ind}{\mathrm{ind}}
\newcommand{\indkg}{\mathrm{ind}^{G}_{KZ}}

\newcommand{\gp}{\mathrm{Gal}\left(\bar{\mathbb{Q}}_p/\mathbb{Q}_p\right)}
\newcommand{\gl}{\mathrm{GL}_{2} (\mathbb{Q}_{p})}
\newcommand{\glfp}{\mathrm{GL}_{2} (\mathbb{F}_{p})}

\newcommand{\symzp}{\mathrm{Sym}^r(\bar{\mathbb{Z}}_{p}^2)}
\newcommand{\symqp}{\mathrm{Sym}^r(\bar{\mathbb{Q}}_{p}^2)}
\newcommand{\vkap}{\bar{V}_{k,a_p}}
\newcommand{\vkapL}{\bar{V}_{k',a_p}}

\newcommand{\fpbar}{\bar{\mathbb{F}}_p}

\begin{document}

\title{Determination of certain mod $p$ Galois representations using local constancy}
\author{Abhik {Ganguli}*}
\author{Suneel {Kumar}**}

\thanks{\noindent Email: *aganguli@iisermohali.ac.in\\
Department of Mathematical Sciences, Indian Institute of Science Education and Research (IISER) Mohali, Sector 81, SAS Nagar,  Punjab-140306, India.}
\thanks{\noindent Email: **suneelm145@gmail.com\\ Department of Mathematical Sciences, Indian Institute of Science Education and Research (IISER) Tirupati, Srinivasapuram, Yerpedu Mandal, Tirupati Dist, Andhra Pradesh-517619, India.}

\begin{abstract}
Let $p \geq 5$ be a prime. Let $k = b + c(p-1)$ be an integer in $[2p+2, p^2 - p +3]$, where $b \in [2,p]$ and $c \in [2, p-1]$. We prove local constancy in the weight space of the mod $p$ reduction of certain two dimensional crystalline representations of $\gp$, where the slope $\nu(a_p)$ is constrained to be in $(1, c)$ and non-integral. We use the mod $p$ local Langlands correspondence for $\text{GL}_{2} (\mathbb{Q}_{p})$ to compute the mod $p$ reductions explicitly, thereby also giving a lower bound on the radius of constancy around the weights $k$ in the above range and under additional conditions on the slope. As an application of local constancy, we obtain explicit mod $p$ reductions at many new values of $k$ and $a_p$.     
\end{abstract}

\maketitle

\section{Introduction }
Let $p$ be an odd prime and $\nu : \bar{\mathbb{Q}}_p^*\rightarrow\mathbb{Q}$ be the normalized $p$-adic valuation such that $\nu(p) = 1$. For an integer $k\geq 2$ and $0\not=a_p\in\bar{\mathbb{Q}}_p$ with slope $\nu(a_p)>0$, let $D_{k,a_p}$ be the weakly admissible filtered $\phi$-module given in \cite{BLZ} with the characteristic polynomial of semilinear Frobenius $\phi$ given by $X^2-a_pX+p^{k-1}$ with jumps in the filtration at $0$ and $k-1$. By the theorem of Colmez-Fontaine (Theorem A in \cite{cf}) there exists a unique irreducible $2$-dimensional crystalline representation $V_{k,a_p}$ of $\gp$ with Hodge-Tate weights $(0,k-1)$ such that $D_{cris}(V_{k,a_p}^*)\cong D_{k,a_p}$. Here $D_{cris}$ is Fontaine's functor defined in \cite{fontainey} and $V_{k,a_p}^*$ is the dual representation of $V_{k,a_p}$. Let $\vkap$ be the reduction of a $\gp$-stable lattice of $V_{k,a_p}$ up to semisimplification. The problem of explicit computation of the mod $p$ reduction $\vkap$ is quite intricate, and a substantial work has been done using local techniques that involve $p$-adic Hodge theory and more recently the mod $p$ local Langlands correspondence for $\text{GL}_{2} (\mathbb{Q}_{p})$ due to Breuil and Berger (\cite{Br03a}, \cite{Br03b}, \cite{BB10}, \cite{B10}). 



We see that ${V}_{k,a_p}$ is completely determined by $a_p$ and the weight $k$ and thus, so is $\vkap$. In this article, fixing $a_p$ we study the question of local constancy of $\vkap$ as a function of $k$ in the weight space. From Berger's local constancy theorem (Theorem B, \cite{Berger12}) we expect local constancy to hold if $k$ and $k'$ are $p$-adically close enough and are in the same class modulo $p-1$. Let $m(k,a_p)$ be the smallest integer $m$ such that $\bar{V}_{k' , a_{p}} \cong \bar{V}_{k , a_{p}}$ for all $k' \in k + p^{m-1}(p-1) \mathbb{Z}_{\geq 0}$.  


We write $k$ as $b + c(p-1) +2 $, where we take $b \in [2,p]$ and $c \geq 0$. The first result giving an explicit upper bound for Berger's constant $m(k,a_p)$ is in Bhattacharya \cite{maam}, and \cite{GK} extends the result significantly to cover more values of $k$ and allowing higher slope. Thus, \cite{maam} \& \cite{GK} give an explicit lower bound on the radius of local constancy. In both \cite{maam} \& \cite{GK} the slope $\nu (a_p) > c$, and the precise reduction $\vkapL$ is also given in the disk of constancy.

In this article, we consider the problem of local constancy in the situation when $\nu(a_p)<c$ and also non-integral. The approach in \cite{maam}, \cite{GK} and our results here is to compute explicitly $\vkapL$ for all $k'$ in the punctured disk $k+p^t(p-1)\mathbb{Z}_{>0}$ for $t > t_{0}$, where $t_{0}$ is given explicitly and sufficiently large. In order to compute $\vkapL$, we use the $\bmod\ p$ local Langlands correspondence for $\gl$. In \cite{maam} and \cite{GK}, Theorem $1.1.1$, \cite{BL} together with Berger's local constancy theorem (Theorem B, \cite{Berger12}) are applied to further determine the reduction $\vkap$ at the center to finally establish local constancy in the disc around the weight $k$. In our present situation, the condition for Berger's local constancy is already satisfied since $\nu(a_p)<c$, proving the existence of local constancy for these small slopes. We use Berger's theorem to infer that since the reduction $\vkapL$ computed in a sufficiently small punctured disk is the same as $\vkap$ at the center, local constancy must hold in the disk $k+p^t(p-1)\mathbb{Z}_{\geq 0}$ for $t > t_{0}$ around the weight $k$.

In this small slope range the lower bound $\nu(a_p) > \lfloor\frac{k-1}{p}\rfloor$ in Bergdall-Levin, \cite{BL} (and also the larger bound from Berger-Li-Zhu, \cite{BLZ}) is not satisfied. Therefore, unlike the generic situation in \cite{maam}, \cite{GK}, we are not able to compute in general the reduction at the center separately, making it difficult to predict the precise reduction in a punctured disk in the small slope regime. Furthermore, as an important application of local constancy in the small slope regime, we instead deduce the precise reduction $\vkap$ at the center in previously unknown cases of weights and slopes (see Corollary \ref{new_reductions}). Let $\nu$ denote $\lfloor\nu(a_p)\rfloor$. Our first result is as follows (Theorem \ref{main_result III} $(1)$):
\begin{theorem}\label{intro 1}
Let $k = b+c(p-1)+2$ with $2\leq c\leq p-1$, $2\leq b\leq p$ and $p\geq 5$. Fix $a_p$ such that $\nu(a_p)$ is non-integral, $1<\nu(a_p)<c-\epsilon$, where $\epsilon\in\{0,1,2\}$ as defined in \eqref{dfn epsilon}, and let $t>\nu(a_p)+c$. If $b\geq c+\nu-1$ such that $b\not=2\nu+1$ and $(b,\nu)\not=(p,1)$, then $\vkapL\cong\ind\left(\omega^{b+1 + \nu (p-1)}_2\right)$ for all $k'\in k+p^t(p-1)\mathbb{Z}_{\geq 0}$. Moreover, the Berger's constant $m(k,a_p)$ exists and $m(k,a_p)\leq \lceil\nu(a_p)\rceil+c+1$. 

\end{theorem}

We note that in Theorem \ref{intro 1} as well as in Theorem \ref{intro 2} below, we treat only the cases where $\vkapL$ is necessarily irreducible. The omitted values of $b$ in the theorems are precisely the possibly reducible cases that arise from the mod $p$ local Langlands correspondence (see Lemma \ref{vrc 1}). We remark that the above theorem can be proved for $1<\nu(a_p)<c+1$ with different methods (see Theorem $1.0.6$ \cite{thesis}). We also observe that the above theorem shows that $\vkapL$ also depends on $\nu$ when $\nu(a_p)<c-\epsilon$ unlike in the cases of $\nu(a_p)>c$ known so far, where $\vkapL$ depends only on the weight $k$ (see Theorem $1.1$ in \cite{GK}). In Theorem \ref{intro 1}, although the condition $b\geq c+\nu-1$ is technical, it arises quite naturally from the broader constraints coming from our method (see \S \ref{proof_strategy} below). 



We next address the question of local constancy in the case when $b\leq c+\nu-2$, wherein we assume Conjecture \ref{mono II conj}. Our next result is as follows (Theorem \ref{main_result III} $(2)$-$(4)$):

\begin{theorem}\label{intro 2}
Let $k = b+c(p-1)+2$ with $2\leq c\leq p-1$, $2\leq b\leq p$ and $p\geq 5$. Fix $a_p$ such that $\nu(a_p)$ is non-integral, $1<\nu(a_p)<c-\epsilon$, and let $t>\nu(a_p)+c$. Assume \textbf{Conjecture \ref{mono II conj}} is true. 
\vspace{0.5em}\\
$(1)$ If $c-1\leq b\leq 2c-4$ and $\nu = c-2$, then $\vkapL\cong\ind\left(\omega_2^{b+1+(b-c+1)(p-1)}\right)$ for all $k'\in k+p^t(p-1)\mathbb{Z}_{\geq 0}$.
\vspace{0.5em}\\
$(2)$   Suppose $2\leq b\leq c-2$ and $1\leq \nu\leq \text{min}\{c-2,p+b-c\}$. Also assume that $b\not=2\nu+2-p$ if $b\leq 2c-2-p$. Then for all $k'\in k+p^t(p-1)\mathbb{Z}_{\geq 0}$, $\vkapL\cong\ind\left(\omega_2^{b+1+(\nu+1)(p-1)}\right)$ if $b\leq \nu$, and $\vkapL\cong\ind\left(\omega_2^{b+1}\right)$ if $b=\nu+1$.
\vspace{0.5em}\\
$(3)$   If $2\leq b\leq 2c-5-p$ and $\nu\in\{p+b-c+1, c-3\}$, then $\vkapL\cong\ind\left(\omega_2^{k_0}\right)$ for all $k'\in k+p^t(p-1)\mathbb{Z}_{\geq 0}$, where $k_0=b+1+(p+b-c+2)(p-1)$.
\vspace{0.5em}\\
In the above cases the Berger's constant $m(k,a_p)$ exists and $m(k,a_p)\leq \lceil\nu(a_p)\rceil+c+1$. 
\end{theorem}


In determining $\vkapL$ in the theorems above, the critical component of our approach is to show certain monomials belong to the kernel of a map $P$ obtained from the $\bmod\ p$ local Langlands machinery (see \S \ref{proof_strategy}). This involves computations with the Hecke operator $T$ requiring a delicate choice of functions, and also some complicated binomial identities. Conjecture \ref{mono II conj} is critical to showing that these monomials (the $q(c)$ in \S \ref{proof_strategy}) are essentially in the kernel of this map $P$. We make an important remark that Theorem \ref{intro 1} is unconditional on Conjecture \ref{mono II conj}. Theorem \ref{proof_conj} proves a substantial portion of Conjecture \ref{mono II conj}. We refer to \S \ref{proof_strategy} for a discussion on this conjecture. 

For slope $0<\nu(a_p)<2$, the $\bmod\ p$ reduction $\vkap$ is completely known for all the weights (see \cite{KBG},\cite{KBGII},\cite{GG15},\cite{BG},
\cite{BGR18},\cite{grai}). Our contribution in computing new cases of $\vkap$ is therefore when the slope $\nu(a_p)>2$. Recall that $b$ is the unique integer in $[2, p]$ such that $k-2 \equiv b \bmod{p-1}$ and $c := \frac{k-2-b}{p-1}$. As a consequence of Theorems \ref{intro 1} \& \ref{intro 2} we have the following corollary that gives the reduction $\vkap$ at new values of $k$ and $a_p$. We note that the first reduction given in each of the three ranges of slopes below are coming from Theorem \ref{intro 1}, and hence unconditional on Conjecture \ref{mono II conj}. 


\begin{cor}\label{new_reductions}
Fix $a_p$ such that $\nu(a_p)$ is non-integral and $1<\nu(a_p)<p-1$. Let $p\geq 13$  and assume \textbf{Conjecture \ref{mono II conj}} is true. Suppose that $k\not\in\{2\nu+3+c(p-1),2\nu+4-p+c(p-1)\}$, where $c$ is given below. 
\begin{enumerate}
    \item If $\nu =2$, then 
    \begin{align*}
        \vkap\cong\begin{cases}
           \ind\left(\omega^{b+1+\nu(p-1)}_2\right)\ &\text{if}\quad k\in\bigcup\limits_{c= 3+\epsilon}^{p-1}[c+3+c(p-1),p+2+c(p-1)]\\ 
           \ind\left(\omega_2^{b+p}\right)\ &\text{if}\quad k = 4p+2\\
            \ind\left(\omega_2^{b+1}\right)\ &\text{if}\quad \{5+c(p-1)|\ 4\leq c\leq p-1\}\\
           \ind\left(\omega_2^{b+1+(\nu+1)(p-1)}\right)\ &\text{if}\quad k\in\{4+c(p-1)|\ 4\leq c\leq p-1\}.
        \end{cases}
    \end{align*}
    \item If $\nu =3$, then 
    \begin{align*}
      \vkap\cong\begin{cases}
         \ind\left(\omega^{b+1+\nu(p-1)}_2\right)\ &\text{if}\quad k\in\bigcup\limits_{c= 4+\epsilon}^{p-1}[c+4+c(p-1),p+2+c(p-1)]\\ 
         \ind\left(\omega_2^{b+1+(b-c+1)(p-1)}\right)\ &\text{if}\quad k\in\{5p+2, 5p+3\}\\
         \ind\left(\omega_2^{b+1}\right)\ &\text{if}\quad k\in\{6+c(p-1)|\ 5\leq c\leq p-1\}\\
         \ind\left(\omega_2^{b+1+(\nu+1)(p-1)}\right)\ &\text{if}\quad k\in\bigcup\limits_{c=5}^{p-1}[4+c(p-1),5+c(p-1)].
      \end{cases}
    \end{align*}
    \item If $4\leq \nu\leq \frac{p-1}{2}$, then
    \begin{align*}
      \vkap\cong\begin{cases}
         \ind\left(\omega^{b+1+\nu(p-1)}_2\right)\ &\text{if}\quad k\in\bigcup\limits_{c= \nu+1+\epsilon}^{p+1-\nu}[c+\nu+1+c(p-1),p+2+c(p-1)]\\
         \ind\left(\omega_2^{b+1+(b-c+1)(p-1)}\right)\ &\text{if}\quad k\in[2+(\nu+2)p, \nu+(\nu+2)p]\\
         \ind\left(\omega_2^{b+1}\right)\ &\text{if}\quad k\in\{\nu+3+c(p-1)|\ \nu+2\leq c\leq p-1\}\\
         \ind\left(\omega_2^{b+1+(\nu+1)(p-1)}\right)\ &\text{if}\quad k\in I,
      \end{cases} 
    \end{align*}
    where $I = \bigcup\limits_{c=\nu+2}^{p-\nu+2}[4+c(p-1),\nu+2+c(p-1)]\bigcup\limits_{c = p-\nu+3}^{p-1}[c+\nu+1-p+c(p-1), \nu+2+c(p-1)]$.
\end{enumerate}
\end{cor}

For a given range of slope, the gaps in the intervals for the weight $k$ in the above corollary are precisely the weights not covered in Theorems \ref{intro 1} \& \ref{intro 2}. Corollary \ref{new_reductions} can be extended to slope $\nu \leq p-3$ with $k$ as given in Corollary \ref{new_k}. We have taken the prime $p \geq 13$ in order to give a more uniform and simplified statement.

We discuss the overlap of our theorems above with known results computing $\vkap$. Firstly, Theorems \ref{intro 1} \& \ref{intro 2} have no overlap with the $(k, a_p)$ covered in \cite{Br03b}, \cite{BLZ} \& \cite{BL}. Theorem \ref{intro 1} has a significant number of cases when the slope is in $(1,2)$, and the $\vkap$ from Theorem \ref{intro 1} above match with those from Theorem $1.1$ in \cite{BG}.
Theorem \ref{intro 2} has an overlap with Theorem $1.1$ in \cite{BG} when $b=2$ (since $\nu(a_p) \in (1,2)$), and the reductions match in this case. Theorem \ref{intro 2} also overlaps with Corollary $1.12$ in \cite{GhR}, precisely when $b \in [2, c-2]$ and $\nu (a_p) \in (b-1, b)\cup  (b, b+1)$ or when $b = c-1$ and $\nu(a_p) = (b-1, b)$. The reductions from Theorem \ref{intro 2} match with the ones coming from \cite{GhR}.  

We next discuss results that apply local methods to compute $\vkap$. Breuil in \cite{Br03b} computes $\vkap$ for weights up to $2p+1$ and for all $a_p$ (see Theorem $3.2.1$ in \cite{B10}). Berger-Li-Zhu in \cite{BLZ} and Bergdall-Levin in \cite{BL} treat this problem for sufficiently large slopes $\nu(a_p)>\lfloor\frac{k-2}{p-1}\rfloor$ and $\nu(a_p)>\lfloor\frac{k-1}{p}\rfloor$ respectively. Buzzard-Gee in \cite{KBG} and \cite{KBGII} determine $\vkap$ for slope in $(0,1)$ and for all weights. 
For slope in $(1,2)$, Ganguli-Ghate in \cite{GG15} compute $\vkap$ up to weights $p^{2}-p$, and Bhattacharya-Ghate in \cite{BG} give the reduction for all weights and slope in $(1,2)$ with an additional assumption at $\nu(a_p)=\frac{3}{2}$. Bhattacharya-Ghate-Rozensztajn in \cite{BGR18} treat the case $\nu(a_p) = 1$ for all weights. We refer to the work of Ghate-Rai in \cite{grai} which treats the case of slope $\nu(a_p)=\frac{3}{2}$. Furthermore, Rozensztajn in \cite{sandra} gives an algorithm to compute $\vkap$ which is efficient for small slopes and weights. We also refer to the recent work of Arsovski \cite{A21} in connection to the slope conjecture of Breuil, Buzzard and Emerton. 
The zig-zag conjecture of Ghate in \cite{ghate} (see also Conjecture $1.1$ in \cite{ghate-zigzag}) gives an explicit description of $\vkapL$ for all half-integral slopes $1/2 \leq \nu(a_p)\leq\frac{p-1}{2}$ and $k' \geq k$ (and $k'$ sufficiently close to $k$ $p$-adically), where $3\leq k \leq p+1$ is such that $k = 2 \nu(a_p) +2 $. This conjecture provides important counterexamples to local constancy in the weight space for certain $a_p$ and $k$ (see Theorem $2.2$ in \cite{ghate}). We remark that in fact $k > 2 \nu(a_p) +2$ in our theorems above. 
We refer to the work of Chitrao-Ghate-Yasuda in \cite{cgy} which explains to a large extent the reason behind these important counterexamples to local constancy. 

For results with alternate approaches using global methods which assume the modularity of the Galois representations, we refer to the work of Deligne, Deligne-Serre and Fontaine-Edixhoven (see [Edi92] and also Theorem $1.2$, \cite{KBG}) for modular forms, and the work of Ganguli \cite{G17} and Bhattacharya-Ganguli \cite{SA} for certain Hilbert modular forms of small slope.


\subsection{Outline of the proof}\label{proof_strategy}
Breuil has constructed a locally algebraic representation $\Pi_{k',a_p}$ of $\gl$ and a $\gl$-stable lattice $\Theta_{k',a_p}$ in $\Pi_{k',a_p}$ such that $\bar{\Theta}_{k',a_p}^{ss}\cong LL(\vkapL)$, where $LL$ is the $\bmod\ p$ local Langlands correspondence (see \S\ref{sec mod p LLC} for more details). Let $r = k'-2\geq 0$ be a non-negative integer. Using the definition of $\Theta_{k',a_p}$, we get a canonical surjection
$P: \indkg(V_r)\twoheadrightarrow\bar{\Theta}_{k',a_p}$ (see \S \ref{sec_filtaration} for more details). For $m\in\mathbb{Z}_{\geq 1}$, we denote 
$V_r^{(m)}:=\{f\in V_r|\ \theta^m \ \text{divides}\ f\ \text{in}\ \fpbar[x,y]\}$,
where $\theta:= x^py-xy^p$. From Remark $4.4$ of Buzzard-Gee in \cite{KBG}, we deduce that $$P :\indkg\left(\frac{V_r}{V_r^{(\nu+1)}}\right)\twoheadrightarrow\bar{\Theta}_{k',a_p}$$ as $\nu(a_p)<\nu+1$, where $\nu = \lfloor \nu (a_p) \rfloor$. We consider the following filtration
\begin{eqnarray*}
 0\subseteq\indkg\left(\frac{V^{(\nu)}_r}{V^{(\nu +1)}_r}\right)\subseteq \indkg\left(\frac{V^{(\nu-1)}_r}{V^{(\nu+1)}_r}\right)\subseteq ...\subseteq\indkg\left(\frac{V_r}{V^{(\nu+1)}_r}\right).
\end{eqnarray*}
For $0\leq m\leq \nu$, observe that $\indkg\left(\frac{V^{(m)}_r}{V^{(m+1)}_r}\right)$ are the successive quotients in the above filtration. By Lemma \ref{lm3.2} the successive quotients $\indkg\left(\frac{V^{(m)}_r}{V^{(m+1)}_r}\right)$ are generated by $[g, F_m(x,y)]$, where $F_m(x,y)=x^{r-m}y^m-x^{r-(s-m)}y^{s-m}$ and $s=k-2$. By a contributing factor, we mean the successive quotient $\indkg\left(\frac{V_r^{(n)}}{V_r^{(n+1)}}\right)$ such that $P:\indkg\left(\frac{V_r^{(n)}}{V_r^{(n+1)}}\right)\twoheadrightarrow\Bar{\Theta}_{k',a_p}$ (we prove that there is always only one contributing factor for Theorem \ref{intro 1} and \ref{intro 2}).

To determine $\bar{\Theta}_{k',a_p}^{ss}$ when $\vkapL$ is irreducible, it is enough to determine the contributing factor (see Lemma \ref{vrc 1}). At first, Proposition \ref{rm_JH_conj} shows that the map $P$ surjects via $\indkg\left(\frac{V_r^{(n_1+1)}}{V_r^{(n_2)}}\right)$, where $n_{1}$ \& $n_{2}$ are as given in the proposition.
Next in Theorem \ref{contr_factor}, we prove that $\indkg\left(\frac{V_r^{(n_0)}}{V_r^{(n_0+1)}}\right)$ is the contributing factor, where $n_0$ is defined in the same theorem. The hypotheses of Theorem \ref{contr_factor} is such that $n_1 = n_0-1$ and $n_2 = n_0+1$, thus giving that contributing factor comes from $n_0$. In particular, the condition $b \geq c+ \nu -1$ in Theorem \ref{intro 1} arises naturally to give us that $n_0 = \nu$. 

We next apply Lemma \ref{vrc 1} together with Theorem \ref{contr_factor} to get the reduction $\vkapL$ (in a puncture disk centered at $k$) in Proposition \ref{vk'ap conj}. Finally we apply Theorem \ref{Berger} in Theorem \ref{main_result III} to determine $\vkap$ at the center, thereby giving us local constancy in the whole disk.   

We now explain our strategy to prove Proposition \ref{rm_JH_conj}.
We do this by showing that $\indkg\left(\frac{V^{(m)}_r}{V^{(m+1)}_r}\right)$ do not contribute to $\bar{\Theta}_{k',a_p}$ for all $m\in[0,n_1]\cup[n_2,\nu]$. We apply Lemma \ref{lm m<c} \& \ref{lm m>c} to the intervals $[0,n_1]$ \& $[n_2,\nu]$ respectively, whenever they are non-empty. To apply Lemma \ref{lm m>c}, we are required to show that $F_m(x,y)\in\kerp$ for all $m\in [n_2,\nu]$. To obtain this, it suffices to prove that $F_m(x,y)\in(V_r^{(m+1)}+\kerp)$ for all $m\in[n_2,\nu]$. This can be seen by reverse induction on $m$ after observing that $F_{\nu}(x,y)\in\kerp$ because $F_{\nu}(x,y)\in(V_r^{(m+1)}+\kerp)$ and $\indkg(V_r^{(\nu+1)})\subset\kerp$.

Next to apply Lemma \ref{lm m<c}, we are required to construct $W_m$ for all $m\in [0,n_1]$ of this lemma. In fact, for $m\geq 1$ it suffices to show $F_m(x,y)\in(V_r^{(m+1)}+\kerp)$ (for $m = 0$ we construct $W_0$ separately). This can be seen by observing that $F_m(x,y)\in(V_r^{(m+1)}+\kerp)$ implies that there exists $v_{m+1}\in V_r^{(m+1)}$ such that $F_m(x,y)-v_{m+1}\in\kerp$. Let $W_m$ be the submodule of $V_r^{(m)}$ generated by $F_m(x,y)-v_{m+1}$. By using Lemma \ref{lm3.2}, and given that $F_m(x,y)-v_{m+1}\in\kerp$, we observe that $W_m$ satisfies the required conditions of Lemma \ref{lm m<c}.

From the above discussion, we deduce that it is enough to prove $F_m(x,y)\in(V_r^{(m+1)}+\kerp)$ for all $m\in[1,n_1]\cup [n_2,\nu]$. Next, we explain the steps to prove $F_m(x,y)\in(V_r^{(m+1)}+\kerp)$ for all $m\in[1,n_1]\cup [n_2,\nu]$. By using Remark $4.4$ in \cite{KBG}, observe that for each $m\leq \nu$ $F_m(x,y)\in(V_r^{(m+1)}+\kerp)$ if and only if $q(c) = x^{r-s+m} y^{s-m} \in(V_r^{(m+1)}+\kerp)$ (as $\nu<\nu(a_p)$; see \S \ref{notations} for the definition of $q(i)$ in general).

Applying Proposition \ref{gen 1}, we deduce that there exists $a_{j,l}\in\mathbb{Z}_p$ such that
$\underset{j}\sum a_{j,l}q(j)\in\kerp$, where $0\leq j\leq c-1 \ \& \ b-m+j(p-1)>0$. Further in Lemma \ref{Vrm cong}, we prove that for all $\epsilon_1\leq j\leq c-m-1$, $q(j)$ are integral linear combinations of $q(j)$ for $c-m\leq j\leq c$. Using this we have $$\underset{1\leq i\leq m+1}\sum \frac{\alpha(i,l)}{p^{\sigma_1(l)}}q(c-m-1+i)\in\kerp,$$ where $\alpha(i,l)$ are as given in \eqref{dfn alpha} and $\sigma_1(l)$ in \eqref{sigma_1l}. This is done in Proposition \ref{gen m 1}. We note from the definition of $\sigma_1{(l)}$ and $\alpha(i,l)$ that $\alpha(i,l)/p^{\sigma_1{(l)}}$ is integral for all $i,l$. Therefore, in order to prove that $q(c)\in(V_r^{(m+1)}+\kerp)$ for all $m\in[1,n_1]\cup [n_2,\nu]$, we need to show that for each $m$ there exists some $l$ such that the following holds:
\begin{enumerate}
    \item The sum $\underset{1\leq i\leq m}\sum \frac{\alpha(i,l)}{p^{\sigma_1(l)}}q(c-m-1+i)$ vanishes modulo $(V_r^{(m+1)}+\kerp)$.
    \item The coefficient $\frac{\alpha(m+1,l)}{p^{\sigma_1(l)}}$ of $q(c)$ (i.e., the last term $i = m+1$) in the above summation is non-zero $\bmod\ p$.
\end{enumerate}
To prove the statement $(1)$ above, we first check if $\frac{\alpha(i,l)}{p^{\sigma_1(l)}}$ vanish modulo $p$ for all $1\leq i\leq m$, but we get a negative answer to this in general. Therefore, we ask whether the monomials $q(c-m-1+i)\in(V_r^{(m+1)}+\kerp)$ for all $1\leq i\leq m$. To answer this question, we consider the matrix $A = (\alpha(i,l))_{\substack{1\leq i\leq m+1\\ 0\leq l\leq m}}$ and check whether the linear systems $AX = e_i(\bmod\ p)$ has a solution for all $1\leq i\leq m$. We note that $\alpha(i,l)$ has a complicated formula, and so proving that the linear system of equations $AX=e_i(\bmod\ p)$ has a solution in general turns out to be quite hard. However, we make the conjecture that the above linear system of equations has a solution (see Conjecture \ref{mat_conj}). We apply Proposition \ref{gen m 1} to translate Conjecture \ref{mat_conj} into Conjecture \ref{mono II conj} which claims that the monomials $q(c-m-1+i)\in(V_r^{(m+1)}+\kerp)$ for all $1\leq i\leq m$. 

Next, we consider the expression $\frac{\alpha(m+1,l)}{p^{\sigma_1(l)}}$ occurring in the statement $(2)$ above. One can indeed check that $\frac{\alpha(m+1,l)}{p^{\sigma_1(l)}} = d_l$, where $d_l$ is as defined in \eqref{dfn_d_l}. Crucially in Lemma \ref{d_l_non_zero}, we prove that for each $m$ given in the lemma, there exists at least one $l$ such that $d_l$ is non-zero modulo $p$. We use this in Proposition \ref{Fm_kerp_conj} to recover the statement $(2)$ above. 
Lemma \ref{d_l_non_zero} is not covering all the values of $m$ in $[1,\nu]$. This is because for the remaining values of $m$, the coefficients $d_l$ are zero modulo $p$ for all $l$. This results in the gap between $n_1$ and $n_2$ as seen in the various subcases of Proposition \ref{rm_JH_conj}, wherein $n_2 > n_1 +1$ in general. In our approach, the vanishing of $d_l$ is the key reason for the restrictions appearing on $b,c$ and $\nu$ in Theorems \ref{contr_factor} \& \ref{main_result III}, the other reason being to force $n_1$ and $n_2$ to be consecutive. The condition that $\nu (a_p)$ is non-integral appears first in Proposition \ref{Fm_kerp_conj}, and is required so that one is able to apply Proposition \ref{gen 1} and Lemma \ref{d_l_non_zero}. 

We note importantly that Proposition \ref{Fm_kerp_conj} requires Conjecture \ref{mono II conj} as a hypothesis, and thus Conjecture \ref{mono II conj} is also needed crucially for Proposition \ref{rm_JH_conj}.
We also remark that Lemma \ref{proof_conj_case1} proves Conjecture \ref{mat_conj} in the cases required to make Theorem \ref{intro 1} (i.e., Theorem \ref{main_result III} $(1)$) unconditional on Conjecture \ref{mono II conj} (see Remark \ref{rk_uncodi_1} for more details).

Finally, we discuss the evidence for Conjecture \ref{mat_conj}. In Theorem \ref{proof_conj}, we prove a substantial portion of Conjecture \ref{mat_conj}. In doing so the binomial identities in Lemmas \ref{mat cong} \& \ref{alpha' cong II} play an important role by simplifying the expression $\alpha(i,l)$. In the cases proved in Theorem \ref{proof_conj}, the matrix $A$ turns out to be non-invertible. In the remaining cases not covered by Theorem \ref{proof_conj}, we verify the conjecture using SageMath code for all values of $b,c,m$ and primes up to $97$. Further, the SageMath computations reveal that the matrix $A$ in the remaining cases are all invertible.

\section{Background}
\subsection{The mod $p$ local Langlands correspondence}\label{sec mod p LLC}
We begin by recalling some notations and definitions.  We fix an algebraic closure $\bar{\mathbb{Q}}_{p}$ of $\mathbb{Q}_{p}$ with the ring of integers $\bar{\mathbb{Z}}_{p}$ and the residue field $\bar{\mathbb{F}}_{p}$. Let $G_p$ and $G_{p^2}$ be the absolute Galois groups of $\mathbb{Q}_p$ and $\mathbb{Q}_{p^2}$ respectively where $\mathbb{Q}_{p^2}$ is the unique unramified quadratic extension of $\mathbb{Q}_p$.  Let $\omega_1 =\omega$ be the mod $p$ cyclotomic character, and $\omega_2$ be a fixed fundamental character of level $2$. We view $\omega_1$ as a character of $\mathbb{Q}^*_p$ via local class field theory, identifying uniformizers with geometric Frobenii. For $a\in\mathbb{Z}^{\geq 0}$ such that $(p+1)\nmid a$ let $\ind(\omega^a_2)$ denote the unique two dimensional irreducible representation of $G_p$ with determinant $\omega^a$ and whose restriction to inertia is isomorphic to $\omega^a_2\oplus\omega^{ap}_2$.

We denote the group $\gl$ by $G$, its maximal compact subgroup $\mathrm{GL}_2\left(\mathbb{Z}_p\right)$ by $K$ and the center of $G$ by $Z\cong \mathbb{Q}^*_p$.  For $r\geq 0$  let $V_r: = \text{Sym}^r(\bar{\mathbb{F}}_{p}^2)$ be the symmetric power representation of $\mathrm{GL}_2(\mathbb{F}_p)$ of dimension $r+1$.  We can also view $V_r$ as representations of $KZ$ by defining the action of $K$ through the natural surjection $K\twoheadrightarrow \mathrm{GL}_2(\mathbb{F}_p)$, and by letting $p$ act trivially. 
For $0\leq r\leq p-1, \ \lambda\in\bar{\mathbb{F}}_p$ and a smooth character $\eta :\mathbb{Q}^*_p\rightarrow \bar{\mathbb{F}}^*_p$,  the representation
$$\pi(r, \lambda, \eta) := \frac{\indkg(V_r)}{T-\lambda}\otimes(\eta\circ\det)$$
is a smooth admissible representation of $G$ where $\indkg$ denotes compact induction (see \cite{Br03a},  \cite{KBG}).  The operator $T$ (see \S \ref{Hecke}) is the Hecke operator $T_p$ generating the Hecke algebra\linebreak $ \text{End}_{G} ( \indkg (V_{r})) = \bar{\mathbb{F}}_{p} [T_p]$. The irreducible subquotients of these representations give all the irreducible smooth admissible representations of $G$ (\cite{BL94},\cite{BL95},\cite{Br03a}). 
For $\lambda\in\bar{\mathbb{F}}_p^{*}$, let $\mu_\lambda$ be the unramified character of $G_p$ that sends the geometric Frobenius to $\lambda$. Then Breuil 's semisimple $\bmod\ p$ local Langlands correspondence $LL$ (see \cite{Br03b}) is as follows:
\begin{itemize}
\item $\lambda = 0:$ \hspace{5em} $\ind(\omega^{r+1}_2)\otimes\eta\xleftrightarrow{ LL}  \pi(r, 0 , \eta)$\\
\item $\lambda\not=0:$\hspace{5em} $(\mu_\lambda\omega^{r+1}\oplus\mu_{\lambda^{-1}})\otimes\eta\xleftrightarrow{LL} \pi(r,\lambda,\eta)^{ss}\oplus\pi([p-3-r], \lambda^{-1},\omega^{r+1}\eta)^{ss}$ \\
where $\{0, 1, ..., p-2\}\ni [p-3-r]\equiv p-3-r\bmod (p-1)$.\\
\end{itemize}

For integers $k\geq 2$ we define $\Pi_{k,a_p}: = \frac{\indkg(\symqp)}{T-a_p}$ as representations of $G$ where $r=k-2$ and $T$ is the Hecke operator from \S \ref{Hecke}. We consider the $G$-stable lattice $\Theta_{k,a_p}$ in the representation $\Pi_{k,a_p}$ (see  \cite{Br03b}, \cite{BB10}) given by 
$$\Theta_{k,a_p}:= \text{image}\left(\indkg(\symzp)\rightarrow\Pi_{k,a_p}\right)\cong\frac{\indkg(\symzp)}{(T-a_p)\indkg(\symqp)\cap\indkg(\symzp)}.$$
By the compatibility of the $p$-adic and mod $p$ local Langlands correspondences (\cite{Br03b},\cite{B10}, \cite{BB10}) we have 
$$\bar{\Theta}^{ss}_{k,a_p}\cong LL(\bar{V}_{k,a_p})\quad \text{where}\quad\bar{\Theta}_{k,a_p}:=\Theta_{k,a_p}\otimes \bar{\mathbb{F}}_p.$$
Since the $\bmod\ p$ local Langlands correspondence is injective,  to determine $\bar{V}_{k,a_p}$ it is enough to compute $\bar{\Theta}_{k,a_p}^{ss}.$

\subsection{Hecke Operator T}\label{Hecke} We give an explicit definition of the Hecke operator $T = T_p$ below (see \cite{Br03b} for more details). For $m = 0$,  set $I_0 = \{0\}$  and for $m>0$, let $I_m = \{[\lambda_0]+p[\lambda_1]+...+p^m[\lambda_{m-1}]\ |\ \lambda_i\in\mathbb{F}_p\}\subset\mathbb{Z}_p$ where square brackets denote Teichm\"{u}ller representatives.  For $m\geq 1$ there is a truncation map $[\ ]_{m-1}: I_m\rightarrow I_{m-1}$ given by taking the first $m-1$ terms in the $p$-adic expansion above.  For $m = 1$, $[\ ]_{m-1}$ is the zero map.  For $m\geq 0$ and $\lambda\in I_m$,  let
$$g^0_{m, \lambda} = \begin{pmatrix}
p^m & \lambda\\
0 & 1
\end{pmatrix}
\quad\text{and }\quad
g^1_{m, \lambda} = \begin{pmatrix}
1 & 0\\
p\lambda & p^{m+1}
\end{pmatrix}. $$
Then we have
$$ G =\underset{\substack{m\geq 0, \lambda\in I_m\\ i\in\{0,1\}}} \coprod KZ(g^i_{m,\lambda})^{-1}.$$
Let $R$ be a $\mathbb{Z}_p$-algebra and $V = \text{Sym}^rR^2$ be the symmetric power representation of $KZ$, modelled on homogeneous polynomials of degree $r$ in the variables $x$ and $y$ over $R$. For $g\in G, \ v\in V$,  let $[g, \ v]$ be the function defined by: $[g, v](g') = g'g\cdot v$ for all $g'\in KZg^{-1}$ and zero otherwise.  Since an element of $\indkg (V)$ is a $V$-valued function on $G$ that has compact support  modulo $KZ$, one can see that every element of  $\indkg(V)$ can be written as a finite sum of $[g, v]$ with $g= g^0_{m\lambda}$ or $g = g^1_{m, \lambda}$, for some $\lambda\in I_m$ and $v\in V$. Then the action of $T$ on $[g, v]$ can be given explicitly when $g = g^0_{n,\mu}$ with $n\geq 0$ and $\mu\in I$.  Let $v = \sum\limits_{j = 0}^{r}c_jx^{r-j}y^j$,  with $c_j\in R$.  We write $T = T^++T^-$ where
\begin{eqnarray*}
T^+([g^0_{n,\mu},v]) &=& \underset{\lambda\in I_1}\sum \left[g^0_{n+1,\mu+p^n\lambda}, \sum\limits_{j = 0}^{r}p^j\left(\sum\limits_{i= j}^{r}c_i{i\choose j}(-\lambda)^{i-j}\right)x^{r-j}y^j\right]\\
T^-([g^0_{n,\mu}, v]) &=& \left[g^0_{n-1, [\mu]_{n-1}}, \sum\limits_{j = 0}^{r}\left(\sum\limits_{i= j}^{r}p^{r-i}c_i{i\choose j}\left(\frac{\mu-[\mu]_{n-1}}{p^{n-1}}\right)^{i-j}\right)x^{r-j}y^j\right]\quad \text{for}\ n>0\\
T^-([g^0_{n,\mu}, v]) &=& \left[\alpha, \sum\limits_{j =0}^{r}p^{r-j}c_jx^{r-j}y^j\right]\quad\text{for}\ n=0,\ \text{where}\ \alpha: = g^1_{0,0}.
\end{eqnarray*}
\subsection{The filtration}\label{sec_filtaration}
Let $k' = k + p^{t}(p-1)d$ satisfy the hypotheses of Theorem \ref{intro 1} \& \ref{intro 2}. Since $t \geq \lceil 2 \nu(a_p) \rceil + \epsilon$, we have $r = k'-2\geq (\nu +1)(p+1)$, where $\nu = \lfloor\nu(a_p)\rfloor$.  From the definition of $V_r$ and $\bar{\Theta}_{k',a_p}$ it follows that there is a natural surjection
 $$P: \indkg(V_r)\twoheadrightarrow\bar{\Theta}_{k',a_p}.$$
Now, let us consider the Dickson polynomial $\theta:= x^py-xy^p\in V_{p+1}.$ Here we note that $\text{GL}_2(\mathbb{F}_p)$ acts on $\theta$ by the determinant character. For $m\in\mathbb{N}$, let us denote 
$$V^{(m)}_r = \{f\in V_r\ | \ \theta^m \ \text{divides}\ f\ \text{ in}\ \bar{\mathbb{F}}_p[x,y]\}$$
which is a subrepresentation of $V_r$.  By using Remark $4.4$ of \cite{KBG}, one can see that the map $P$ factors through $\indkg\left(\frac{V_r}{V^{(\nu+1)}_r}\right)$, where $\nu: = \lfloor\nu(a_p)\rfloor$.  So let us consider the following chain of submodules
\begin{eqnarray}{\label{filt}}
0\subseteq\indkg\left(\frac{V^{(\nu)}_r}{V^{(\nu +1)}_r}\right)\subseteq \indkg\left(\frac{V^{(\nu-1)}_r}{V^{(\nu+1)}_r}\right)\subseteq ...\subseteq\indkg\left(\frac{V_r}{V^{(\nu+1)}_r}\right).
\end{eqnarray}
 For $0\leq m\leq \nu$,  observe that $\indkg\left(\frac{V^{(m)}_r}{V^{(m+1)}_r}\right)$ are the successive quotients in the above filtration.  In the following two lemmas we make precise the notion of a successive quotient not contributing to $\bar{\Theta}_{k',a_p}$ via the map $P$. 
 
\begin{lemma}{\label{lm m<c}}
Let $1\leq n\leq \nu $ and assume for all $0\leq m\leq n-1$ that there exists $W_m\subset V^{(m)}_r$ such that $P\left(\ind^{G}_{KZ}(W_m)\right) = 0$ and $W_m\twoheadrightarrow \frac{V^{(m)}_r}{V^{(m+1)}_r}$ induced by the inclusion $W_m\subset V^{(m)}_r$. Then the map $P$ restricted to $\ind^{G}_{kZ}\left(\frac{V^{(n)}_r}{V^{(\nu+1)}_r}\right)$ is a surjection.
\end{lemma}

\begin{lemma}{\label{lm m>c}}
Let $1\leq n\leq \nu $ and suppose for all $n\leq m\leq \nu$ that there exists $G_m(x,y)\in V_r$ such that $P([g, \ G_m(x,y)]) = 0$. If $G_m(x,y)$ generates $\frac{V^{(m)}_r}{V^{(m+1)}_r}$ then the map $P$ factors through  $\ind^{G}_{kZ}\left(\frac{V_r}{V^{(n)}_r}\right)$.
\end{lemma}

\subsection{JH factors of ${V^{(n)}_r}/{V^{(n+1)}_r}$}
In this subsection, we assume simply that $r, n \in \mathbb{Z}_{\geq 0}$ such that $r \geq (n+2)(p+1) -3$ so that we can apply Theorems $(4.1)$ and $(4.2)$ of \cite{glover} below. 
Next, we determine the Jordan-Holder factors of the successive quotients $\frac{V^{(n)}_r}{V^{(n+1)}_r}$. Let $D$ denote the determinant character of $\glfp$. Let us write $r-n(p+1) = r'+d'(p-1)$ such that $p\leq r'\leq 2p-2$ and for some $d'\in\mathbb{Z}^{\geq 0}$.  By Theorems $(4.1)$ and $(4.2)$ of \cite{glover} together with Lemma $5.1.3$ of \cite{Br03b} gives:\\
(i) if $r' = p$ then
\begin{eqnarray}{\label{r' = p}}
0\longrightarrow V_1\otimes D^n\longrightarrow\frac{V^{(n)}_r}{V^{(n+1)}_r}\longrightarrow V_{p-2}\otimes D^{n+1}\longrightarrow 0.
\end{eqnarray}
The first map sends $(x, \ y)$ to $(\theta^n x^{r-n(p+1)}, \ \theta^n y^{r-n(p+1)})$ and the second map sends $\theta^n x^{r-n(p+1)-1}y$ to $x^{p-2}$.\\
(ii) if $r'\not = p$ then 
\begin{eqnarray}{\label{r' not p}}
0\longrightarrow V_{r'-(p-1)}\otimes D^n\longrightarrow\frac{V^{(n)}_r}{V^{(n+1)}_r}\longrightarrow V_{2(p-1)-r'}\otimes D^{n+r'-(p-1)}\longrightarrow 0.
\end{eqnarray}
The first map sends $(x^{r'-(p-1)}, \ y^{r'-(p-1)})$ to $(\theta^n x^{r-n(p+1)}, \ \theta^ny^{r-n(p+1)})$ because ${{r'}\choose p-1}\equiv 0\bmod p$ as $1\leq r'-p\leq p-2$. For $r'-(p-1)\leq i\leq p-1$, the second map sends $\theta^nx^{r-n(p+1)-i}y^i$ to $\alpha_i\ x^{p-1-i}y^{p-1-r'+i}$ where $\alpha_i:= (-1)^{r'-i}{{2(p-1)-r'}\choose p-1-r'+i}\not\equiv 0\bmod p$ because $0\leq 2(p-1)-r'\leq p-3$ and $0\leq p-1-r'+i\leq 2(p-1)-r'$.\\

\subsection{Some crucial results} 
In this section, we state Berger's local constancy theorem and some crucial lemmas required later.
\begin{theorem}\textup{[Berger \cite{Berger12}]}{\label{Berger}}
Suppose $a_p\not= 0$ with $\nu (a_{p}) > 0$ and $k>\footnote{\text{Theorem B in \cite{Berger12} is actually stronger, wherein the lower bound on $k$ is $3\nu(a_p)+\alpha(k-1)+1$.}\linebreak \text{Here $\alpha(k-1)=\underset{n\geq 1}\sum \lfloor(k-1)/(p^{n-1}(p-1))\rfloor$. We note that although $\frac{(k-1)p}{(p-1)^2} \geq \alpha(k-1) $, it is easier}\linebreak \text{to use this bigger bound for calculations.}}3\nu(a_p)+\frac{(k-1)p}{(p-1)^2}+1$, then there exists $m = m(k,a_p)$ such that $\vkapL\cong\vkap$, if $k'-k\in p^{m-1}(p-1)\mathbb{Z}_{\geq 0}$. 
\end{theorem}
\vspace*{2em}
For integers $0\leq m\leq s$ let us define polynomials $F_m$ in  $V_r$ as follows
$$F_m(x,y) := x^my^{r-m}-x^{r-s+m}y^{s-m}$$
where $r>s$ and $r\equiv s\bmod (p-1)$.
\begin{lemma}[Bhattacharya, Lemma 3.2, \cite{maam}]{\label{lm3.2}}
Let $r\equiv s\bmod (p-1)$, and $t = \nu(r-s)\geq 1$ and $1\leq m\leq p-1$.\begin{enumerate}
\item For $s\geq 2m$, the polynomial $F_m$ is divisible by $\theta^m$ but not by $\theta^{m+1}$.
\item For $s>2m$, the image of $F_m$ generates the subquotient $\frac{V^{(m)}_r}{V^{(m+1)}_r}$ as a $\glfp$-module.
\end{enumerate}
\end{lemma}
\begin{lemma}[Lemma 6.1, \cite{GK}]{\label{vrc 1}}
Let $p \geq 5$ and $r, n$ be integers such that $0\leq n\leq p-1$ and $r \geq (n+2)(p+1)-3$. Let $b \equiv r \bmod p-1$ such that $2 \leq b \leq p$. Suppose the map
\begin{equation}{\label{map p}}
 P:\ind^G_{KZ}\left( \frac{V^{(n)}_r}{V^{(n+1)}_r}\right)\rightarrow\bar{\Theta}_{r+2,a_p}
\end{equation}
is a surjection. Further if $(b,n)\not\in\{(p-2, \ 0), (p,0),(p,1)\}$ and also $b\not\in\{2n\pm 1, 2(n+1)-p, 2n-p\}$ then 
$$\bar{V}_{r+2,a_p}\cong \begin{cases}
\ind\left(\omega^{b+ n(p-1) +1}_2\right) & \text{if}\quad 2n+1\leq b\leq p\\
\ind\left(\omega^{b+ (n+1)(p-1) +1}_2\right) & \text{if}\quad 2n+1-(p-1)\leq b\leq 2n\\
\ind\left(\omega^{b+ (n+2)(p-1) +1}_2\right) & \text{if}\quad 2(n+1)-2(p-1)\leq b\leq 2n-(p-1).\\
\end{cases}$$
\end{lemma}

\subsection{Notations and Conventions}\label{notations}
We fix the following conventions in the rest of this article unless stated otherwise:
\begin{enumerate}
\item The integer $p$ always denotes a prime number greater than or equal to $5$.  The integers $b$ and $c$ are from $\{2, 3, ...,p\}$ and $\{0,1, ..., p-2\}$ respectively.
\item We define $\epsilon$ as follows 
\begin{eqnarray}{\label{dfn epsilon}}
 \epsilon =\begin{cases} 
 0 & \text{if}\quad 2c-1\leq b\leq p\\
 1 & \text{if} \quad 2(c-1)-p\leq b\leq 2(c-1)\\
 2 & \text{if}\quad 2\leq b\leq 2(c-1)-(p+1).
\end{cases}
\end{eqnarray}
\item  We define $\epsilon_1$ as follows
\begin{align}\label{epsilon II}
\epsilon_1 = \begin{cases}
0\ & \text{if}\ 2m+1\leq b\leq p\\
1 & \text{if}\ 2m+1-(p-1)\leq b\leq 2m\\
2 & \text{if}\ 2\leq b\leq 2m-(p-1).
\end{cases}
\end{align}
\item We write $s = b+c(p-1)$ and $r = s+p^t(p-1)d$ with $p\nmid d$,  and $t,d\in\mathbb{N}$ and so $s < r$.
\item For $n\in\mathbb{Z}_{\geq 0}$ and $k\in\mathbb{Z}$,  we define ${n\choose k} = 0$ if $k > n$ or $k < 0$ and the usual binomial coefficient otherwise.
\item For $A\equiv B$, where $A,B\in \text{M}_n(\mathbb{Z}_p)$ we mean that $A\equiv B\bmod p$.
\item Unless stated otherwise,  for $A, B\in\indkg(\symqp)$, by $A\equiv B$ or $A\equiv B\bmod p$ we mean that $A- B$ is in $\textbf{m}_{\bar{\mathbb{Z}}_p}\indkg(\symzp)$.
\item By the vectors $\{\textbf{e}_j\}$ we mean the standard basis of a free module over $\mathbb{Z}_p$. 
\item For $v\in \text{Sym}^r(\bar{\mathbb{F}}^2_p)$, by $v\in\kerp$ or $v \in V^{m+1}_r+\kerp$ we mean that $[1,\ v]\in\kerp$ or $[1,\ v] \in \indkg (V^{m+1}_r) + \kerp$ respectively.
\item We define $q(i)= x^{r-b+m-i(p-1)}y^{b-m+i(p-1)}$ for all $n_0\leq i\leq c\ \&\ 1\leq m\leq p-1$, where $n_0 = 0$ if $b\geq m$ and $1$ otherwise.
\end{enumerate}

\section{Binomial Identities}

In this section, we prove some technical lemmas, which we use later.
					\begin{lemma}{\label{dfn beta}}
						For $0\not = a,m\in\mathbb{Z}_{\geq 0}$ and $1\leq i\leq m+1$, let 
						$$\beta(a,i):= \begin{cases}
							{{m+1}\choose i} & \text{if}\ a = 1\\
							\underset{1\leq l\leq a-1}\sum(-1)^{l+1}{{m+1}\choose l}\beta(a-l, i)+(-1)^{a-1}{{m+1}\choose i+a-1} &\text{if}\ a\geq 2
						\end{cases}$$
						then $\beta(a,i) = {{i+a-2}\choose a-1}{{m+a}\choose i+a-1}$.            
					\end{lemma}

					\begin{proof}
						We prove the above result by induction on $a$. Observe that it is trivially true for $a = 1$ and for all $i$. By induction, assume the above result is true for all $1\leq a\leq k-1$ and we prove it for $a = k\geq 2$. Therefore, we need to prove 
						$${{i+k-2}\choose k-1}{{m+k}\choose i+k-1}=\underset{1\leq l\leq k-1}\sum(-1)^{l+1}{{m+1}\choose l}\beta(k-l, i)+(-1)^{k-1}{{m+1}\choose i+k-1}.$$
						Using induction, we get 
						\begin{align}\label{beta_equivalent}
						    &\underset{0\leq l\leq k-1}\sum (-1)^l{{m+1}\choose l}{{i+k-l-2}\choose i-1}{{m+k-l}\choose m+1-i} = (-1)^{k-1}{{m+1}\choose i+k-1}.
						\end{align}
                        Hence, we need to prove the above equation to prove our lemma.
						Now,
						$$(x-1)^{m+1} = \underset{0\leq l\leq m+1}\sum(-1)^l{{m+1}\choose l}x^{m+1-l}.$$
						Multiply the above equation by $x^{k-1}$ ($k\geq 2$), we get 
						\begin{equation*}
							(x-1)^{m+1}x^{k-1} = \underset{0\leq l\leq m+1}\sum(-1)^l{{m+1}\choose l}x^{m+k-l}.
						\end{equation*}
						By differentiating the above equation with respect to $x$, $(m+1-i)$ times and multiplying by $\frac{x^{-1}}{(m+1-i)!}$,  we get 
						\begin{eqnarray*}
							\underset{0\leq j\leq m+1-i}\sum {{m+1}\choose m+1-i-j}{{k-1}\choose j}(x-1)^{i+j}x^{k-2-j}
						\end{eqnarray*}
						\begin{eqnarray*}
							= \underset{0\leq l\leq m+1}\sum (-1)^l{{m+1}\choose l}{{m+k-l}\choose m+1-i}x^{i+k-l-2}.
						\end{eqnarray*}
						Note that in the first sum, $j$ can be taken over the range $0\leq j\leq n_1$ where $n_1:= \text{min}\{m+1-i,\ k-1\}$. The last sum can be taken over the range $0\leq l\leq i+k-1$, since ${{m+k-l}\choose m+1-i} = 0$ if $i+k-1<l\leq m+1$ and ${{m+1}\choose l} = 0$ if $m+1<l\leq i+k-1$. Thus, we have 
						\begin{eqnarray*}
							\underset{0\leq j\leq n_1}\sum {{m+1}\choose m+1-i-j}{{k-1}\choose j}(x-1)^{i+j}x^{k-2-j} \\
							= \underset{0\leq l\leq i+k-1}\sum (-1)^l{{m+1}\choose l}{{m+k-l}\choose m+1-i}x^{i+k-l-2}.
						\end{eqnarray*}
						Again, by differentiating the above equation with respect to $x$, $(i-1)$ times and dividing by $(i-1)!$, we get
						\begin{align*}
							&	\underset{0\leq j\leq n_1}\sum {{m+1}\choose i+j}{{k-1}\choose j}\left(\underset{0\leq j'\leq i-1}\sum{{i+j}\choose i-1-j'}(x-1)^{j+j'+1}\frac{1}{j'!}d^{j'}(x^{k-2-j})\right)\\
							&= \underset{0\leq l\leq i+k-2}\sum (-1)^l{{m+1}\choose l}{{m+k-l}\choose m+1-i}{{i+k-l-2}\choose i-1}x^{k-l-1}+(-1)^k{{m+1}\choose i+k-1}x^{-i}.
						\end{align*}
						By putting $x=1$ in the above equation, we get 
						$$\underset{0\leq l\leq i+k-2}\sum (-1)^l{{m+1}\choose l}{{m+k-l}\choose m+1-i}{{i+k-l-2}\choose i-1}+(-1)^k{{m+1}\choose i+k-1} = 0$$
						as $(x-1)^{j+j'+1}=0$ at $x=1$ (since $j,j'\geq 0$). Observe that the above summation can be taken over the range $0\leq l\leq k-1$ as ${{i+k-l-2}\choose i-1} = 0$ if $k-1<l\leq i+k-2$. Thus, we obtained the required equation \eqref{beta_equivalent}.
					\end{proof} 

     For $0\leq l\leq p-1$ and $1\leq i\leq m+1$, we define $\alpha(i,l)$ as follows
					\begin{eqnarray}{\label{dfn alpha}}
						\alpha(i,l):= \begin{cases}
							\alpha_1(i,l)+{{r-l}\choose b-m+(i+c-m-1)(p-1)} & \text{if}\ 1\leq i\leq m\\
							\alpha_1(i,l) & \text{if}\  i= m+1
						\end{cases}
					\end{eqnarray}
					where 
					\begin{align}\label{dfn alpha_1}
					    \alpha_1(i,l) =(-1)^{i+1}\underset{1\leq a\leq c-m-\epsilon_1}\sum{{r-l}\choose b-m+(c-m-a)(p-1)}\beta(a,i)
					\end{align}
					and $\beta(a,i)$ is defined as in Lemma \ref{dfn beta}. For $ 1\leq i\leq m+1$, we define $\alpha'(i,l)$ as follows 
					\begin{eqnarray}{\label{dfn alph'}}
				    &\alpha'(i,l):= \begin{cases}
                                         0 \quad \text{if}\ \ 0\leq l\leq b-c\\
                                       (-1)^{i+1}\underset{c_0\leq a\leq c_1}\sum{{p+b-c-l}\choose b-c+a}{{c-1}\choose c-m-a}\beta(a,i)\ \text{if}\ c_2\leq l\leq p+b-c\\
                                      (-1)^{i+1}\underset{1\leq a\leq c_1}\sum{{2p+b-c-l}\choose p+b-c+a}{{c-2}\choose c-m-1-a}\beta(a,i)\ \text{if}\ p+c_2\leq l\leq m
                                      \end{cases}\hspace{2.3em}
					\end{eqnarray}
                    where $c_0=\text{max}\{c-b,1\}$, $c_1=c-m-\epsilon_1$, $c_2=b-c+1$.
                    \begin{lemma}{\label{alpha cong}}
						Let $r= s+p^t(p-1)d$, with $p\nmid d$, $s =b+c(p-1)$ and $t\geq 2$.  Assume that $2\leq b\leq p$, $0\leq c\leq p-1$,  $1\leq m\leq c-1,\ 1\leq i\leq m+1$ and $0\leq l\leq m$. Let $$X(i,l)={{r-l}\choose b-m+(i+c-m-1)(p-1)}.$$
                         Then $\alpha_1(i,l)\equiv \alpha'(i,l)$ and
                         \begin{eqnarray*}
                             &X(i,l)\equiv \begin{cases}
                            {{b-c-l}\choose b-c+1-i}{{c}\choose c-m-1+i} &\text{if}\  \substack{0\leq l\leq b-c,\ 1\leq i\leq b-c+1}\vspace{.35em}\\
                            {{p+b-c-l}\choose b-c+1-i}{{c-1}\choose c-m-1+i}& \text{if}\  \substack{b-c+1\leq l\leq p+b-c,\ 1\leq i\leq b-c+1}\\
                            0&\text{if}\ \substack{0\leq l\leq b-c,\ b-c+2\leq i\leq p+b-c+1}\\
                            {{p+b-c-l}\choose p+b-c+1-i}{{c-1}\choose c-m-2+i}& \text{if}\ \substack{b-c+1\leq l\leq p+b-c,\ b-c+2\leq i\leq p+b-c+1}\vspace{.35em}\\
                            {{2p+b-c-l}\choose p+b-c+1-i}{{c-2}\choose c-m-2+i}& \text{if}\ \substack{p+b-c+1\leq l\leq p-1,\ b-c+2\leq i\leq p+b-c+1}\\
                            0&\text{if}\ \substack{b-c+1\leq l\leq p+b-c,\ p+b-c+2\leq i\leq p-1}\\
                            {{2p+b-c-l}\choose 2p+b-c+1-i}{{c-2}\choose c-m-3+i}& \text{if}\ \substack{p+b-c+1\leq l\leq p-1,\ p+b-c+2\leq i\leq p-1.}\\
                         \end{cases}
                         \end{eqnarray*}
                        \end{lemma}
					\begin{proof}
						First, we note that 
						\begin{eqnarray*}
							r-l & = &b-c-l+cp+p^t(p-1)d\\
							b-m+(c-m-a)(p-1) &=& b-c+a+(c-m-a)p\\
							b-m+(c-m-1+i)(p-1) &=& b-c+1-i+(c-m-1+i)p
						\end{eqnarray*}
						If $0\leq l\leq b-c$ and $1\leq i\leq m+1$, then by Lucas Theorem we have 
						\begin{eqnarray*}
							&&(-1)^{i+1}\underset{1\leq a\leq c-m}\sum{{r-l}\choose b-m+(c-m-a)(p-1)}\beta(a,i)\\
							&&\equiv (-1)^{i+1}\underset{1\leq a\leq c-m}\sum{{b-c-l}\choose b-c+a}{c\choose c-m-a}{{i+a-2}\choose a-1}{{m+a}\choose i+a-1}\\
							&&\equiv 0\bmod p.
						\end{eqnarray*}
						The last congruence follows since ${{b-c-l}\choose b-c+a} = 0$ as $b-c-l<b-c+a$. 
						Thus, we obtain the result in this case by using the above calculation together with $${{r-l}\choose b-m+(c-m-a+i)(p-1)}\equiv {{b-c-l}\choose b-c+1-i}{c\choose c-m-1+i}.$$
						By a similar computation, one obtains results in all other cases also.
					\end{proof}

Let $S=\{j\ |\ 0\leq j\leq c-1 \&\ b-m+j(p-1)>0\}$. For $0\leq l\leq p-1$, we define 
\begin{align}\label{sigma_1l}
\nonumber y(j,l)&={{r-l}\choose b-m+j(p-1)}\ \ \forall\ \ j\in S\vspace{1.5em}\\
\&\ \sigma_1(l)& =\begin{cases}
1 & \quad \text{if}\ \ y(j,l)\equiv 0\bmod p\ \forall\ j\in S\\
0 & \quad \text{if}\ \ y(j,l)\not\equiv 0\bmod p\ \text{for at least one}\ j\in S.
\end{cases} 
\end{align}
For any $1\leq \nu\leq p-1$ and $0\leq l\leq \nu-\sigma_1(l)$, we define $d_l$ as follows
\begin{align}\label{dfn_d_l}
d_l:=\frac{(-1)^m}{p^{\sigma_1(l)}}\underset{\epsilon_1\leq j\leq c-m-1}\sum {{c-1-j}\choose m}{{r-l}\choose b-m+j(p-1)}.
\end{align}
			   
\begin{lemma}\label{d_l_non_zero}
Let $s= b+c(p-1)$ with $2\leq b\leq p, \ 1\leq c\leq p-1$ and  $r= s+p^t(p-1)d$, with $t\geq 2$, $p\nmid d$. Fix a $\nu$ in $[1, c-1- \epsilon]$. Assume that $1\leq m\leq \nu $ and $l\in[0,\nu-\sigma_1(l)]$. 
\begin{enumerate}
\item If $c-1\leq b\leq p$, and $1\leq m\leq \text{min}\{\nu-1,b-c\}$, then for each $m$ there exists at least one $l$ such that $d_l\not=0\bmod p$.
\item  If $c-1\leq b\leq p$, and $b-\nu \leq m\leq \nu$, then for each such $m$ there exists at least one $l$ such that $d_l\not=0\bmod p$.
\item If $2c-2-p\leq b\leq c-2$, and $1\leq \nu\leq c-2$, then for each $m\in([b-\nu,\nu]$\linebreak$\cap[1,b-1])\cup[b,\nu-1]$, there exists at least one $l$ such that $d_l\not=0\bmod p$.
\item If $2\leq b\leq 2(c-1)-p-1$, and $1\leq \nu\leq p+b-c$, then for each $m\in([b-\nu,\nu]\linebreak\cap[1,b-1])\cup[b,\nu-1]$, there exists at least one $l$ such that $d_l\not=0\bmod p$. 
\item If $2\leq b\leq 2(c-1)-p-1$ and $p+b-c+1\leq \nu\leq c-3$, then for each $m\in[1,p+b-c]\cup[p+b-\nu-1, \nu]$, there exists at least one $l$ such that $d_l\not=0\bmod p$.
\end{enumerate}
\end{lemma}
\begin{rk}
 The range of $l$ is $[0,\nu-\sigma_1(l)]$ in Lemma \ref{d_l_non_zero} since it is applied in Proposition \ref{Fm_kerp_conj} wherein $l$ is in $[0, \nu]$ (or $[0, \nu-1]$). Further, in almost all cases of above lemma, the $l$ for which $d_l\not\equiv 0\bmod p$ is given in terms of $m$. As will be clear in the proof below, the range of $m$ in the cases below is given primarily so that the desired $l$ (such that $d_l\not\equiv 0\bmod p$) lies in $[0,\nu-\sigma_1(l)]$.     
\end{rk}
\begin{proof}
For each $b,c,m$ as in the statement, in the following cases, we prove that there exists at least one $l\in[0,\nu-\sigma_1(l)]$ such that $d_l\not\equiv 0\bmod p$. For the above range of $0\leq l\leq \nu-\sigma_1(l)$ and $\epsilon_1\leq j\leq c-m-1$, we define  
$$a_{j,l}=\frac{{{r-l}\choose b-m+j(p-1)}}{p^{\sigma_1(l)}}.$$

\textbf{Case (1)} $c-1\leq b\leq p$ and $1\leq m\leq \text{min}\{\nu-1,b-c\}$\\
Note that actually $b \geq c$ in this case. We prove that there exists at least one $l\in[m,\ \text{min}\{\nu-1,b-c\}]$ such that $d_l\not\equiv0(\bmod p)$. By Lemma $3.4$ in \cite{GK}, we have
\begin{align*}
{{r-l}\choose b-m+j(p-1)}&\equiv {{b-c-l}\choose b-m-j}{c\choose j}\quad\text{if}\ 0\leq l\leq b-c\ \&\ 0\leq j\leq c-1\ (\text{as}\ c\leq b-m)\\
&\equiv 0 \quad\text{if}\ m\leq l\leq b-c\ \&\ 0\leq j\leq c-1.
\end{align*}
Therefore, $\sigma_1(l) = 1$ for all $m\leq l\leq\text{min}\{\nu-1,b-c\}$. Using Lemma $3.5$ in \cite{GK}, we get
\begin{align*}
\frac{{{r-l}\choose b-m+j(p-1)}}{p}&\equiv \frac{(-1)^{l-m}{{b-m}\choose j}{{p-1+m-l}\choose c-1-j}}{{{b-m-c}\choose l-m}{{b-m}\choose c}}\quad\ \text{if}\ 0\leq j\leq c-1\ \& \ m\leq l\leq \text{min}\{\nu-1,b-c\}.
\end{align*}

Next, observe that $\epsilon_1 = 0$ as $b\geq m+c\geq 2m+1+\epsilon$. Hence, for $m\leq l\leq \text{min}\{\nu-1,b-c\}$ we have
\begin{align*}
d_l&\equiv (-1)^l\underset{0\leq j\leq c-m-1}\sum \frac{{{b-m}\choose j}{{p-1+m-l}\choose c-1-j}{{c-1-j}\choose m}}{{{b-m-c}\choose l-m}{{b-m}\choose c}}\\
&\equiv \frac{(-1)^l{{p-1+l-m}\choose m}}{{{b-m-c}\choose l-m}{{b-m}\choose c}}\underset{0\leq j\leq c-m-1}\sum {{b-m}\choose j}{{p-1-l}\choose c-m-1-j}.
\end{align*}
Here, the last congruence follows as 
\begin{align}\label{case1_bino_rearr}
    {{p-1+m-l}\choose c-1-j}{{c-1-j}\choose m} = {{p-1+m-l}\choose m}{{p-1-l}\choose c-m-1-j}.
\end{align}
Using Vandermonde's identity, we get
$$d_l\equiv \frac{(-1)^l{{p-1+m-l}\choose m}{{p+b-(m+1+l)}\choose c-m-1}}{{{b-m-c}\choose l-m}{{b-m}\choose c}}.$$
By taking $l = m$, we get 
\begin{align*}
    d_m&\equiv \frac{(-1)^m{{p-1}\choose m}{{p+b-(2m+1)}\choose c-m-1}}{{{b-m}\choose c}}\\
    &\equiv\frac{(-1)^m{{p-1}\choose m}{{b-(2m+1)}\choose c-m-1}}{{{b-m}\choose c}}\quad \text{(as $b\geq 2m+1$).}
\end{align*}
Hence, we get $d_m\not\equiv 0(\bmod\ p)$ as $b-2m-1-(c-m-1)=b-c-m\geq 0$ and $b-2m+1\leq p-1$.\\
\textbf{Case (2)} $c-1\leq b\leq p$ and $ b-\nu \leq m\leq \nu$ (so $\text{max}\{1,b-c+1\}\leq m$)\\
Observe that $b\leq 2c-3$ as for $b\geq 2c-2$ the range of $m$ is empty. 
In this case, we claim that there exists at least one $l\in[b-c+1, \nu]\subset[0,\nu]$ such that $d_l\not\equiv0(\bmod p)$ as $\sigma_1(l) = 0$ for all $l\in[b-c+1, \nu]$. We take $j = b-m$ to show that $\sigma_1(l) = 0$ (also using that $b-c+1 \leq m \leq c-1-\epsilon$) in \eqref{sigma_1l}. 
Next, by Lemma $3.4$ in \cite{GK} we have
$$a_{j,l}\equiv {{p+b-c-l}\choose b-m-j}{{c-1}\choose j}\quad\text{if}\ \ b-c+1\leq l\leq \nu \leq p+b-c\ \&\ \epsilon_1\leq j\leq c-m-1.$$
Hence, we have 
\begin{align*}
d_l &\equiv (-1)^m\underset{\epsilon_1\leq j\leq c-m-1}\sum {{p+b-c-l}\choose b-m-j}{{c-1}\choose j}{{c-1-j}\choose m}.
\end{align*}
Since \begin{align}\label{case2_bino_rearr}
    {{c-1}\choose j}{{c-1-j}\choose m} = {{c-1}\choose m}{{c-m-1}\choose j}
\end{align}
so, we have
\begin{align*}
d_l&\equiv (-1)^m{{c-1}\choose m}\underset{\epsilon_1\leq j\leq c-m-1}\sum {{p+b-c-l}\choose b-m-j}{{c-m-1}\choose j}\\
&\equiv  (-1)^m{{c-1}\choose m}\underset{\epsilon_1\leq j\leq b-m}\sum {{p+b-c-l}\choose b-m-j}{{c-m-1}\choose j}.
\end{align*}
The last congruence follows as $c-m-1\leq b-m$ and ${{c-m-1}\choose j} =0$ for all $c-m-1<j\leq b-m$.
Using Vandermonde's identity, we get
\begin{align*}
d_l\equiv\begin{cases}
(-1)^m{{c-1}\choose m}{{p+b-m-1-l}\choose b-m}\quad &\text{if}\ b\geq 2m+1\\
(-1)^m{{c-1}\choose m}\left({{p+b-m-1-l}\choose b-m}-{{p+b-c-l}\choose b-m}\right)\quad &\text{if}\ b\leq 2m.
	\end{cases}
\end{align*}
If $b\geq 2m+1$, then 
\begin{align*}
d_l&\equiv \begin{cases}
(-1)^m{{c-1}\choose m}{{b-m-1-l}\choose b-m}\quad &\text{if}\ \ b-c+1\leq l\leq b-m-1\\
(-1)^m{{c-1}\choose m}{{p-1}\choose b-m}\quad&\text{if}\ \ l = b-m
\end{cases}
\end{align*}
\begin{align*}
\implies\ \ d_l\equiv \begin{cases}
0\quad &\text{if}\ \ b-c+1\leq l\leq b-m-1\\
(-1)^m{{c-1}\choose m}{{p-1}\choose b-m}\quad&\text{if}\ \ l = b-m.
\end{cases}
\end{align*}
If $b\leq 2m$ and $0\leq l = b-m-1$, then 
\begin{align*}
d_l &\equiv (-1)^m{{c-1}\choose m}\left({{p}\choose b-m}-{{p-(c-m-1)}\choose b-m}\right)\\
&\equiv (-1)^{m+1}{{c-1}\choose m}{{p-(c-m-1)}\choose b-m}\not\equiv 0.
\end{align*}
From the above calculation, we get $d_l\not\equiv0\bmod p$ for $l = b-m-1$ if $b\leq 2m$ and $l = b-m$ if $b\geq 2m+1$. When $m\geq b-\nu$, we note that $l = b-m-1$ and $l = b-m$ lie in the required range $[0,\nu]$. Hence, for each $m\in[b-\nu,\nu]$ there exists at least one $l$ such that $d_l\not\equiv 0\bmod p$.\\
\textbf{Case (3)} $2c-2-p\leq b\leq c-2$ and $1 \leq \nu\leq c-2$\\
Note that in this case $\epsilon_1\in\{0,1\}$ as $b\geq 2(c-2)+1-(p-1)\geq 2m+1-(p-1)$ for all $m\leq c-2$.\\
\textbf{Subcase (I)} $m\in[b-\nu, \nu]\cap[1,b-1]$\\
For $1\leq m\leq b-1$, we observe that $\sigma_1(l) = 0$ for all $l\in[0,\nu]$ (take $j=b-m$ as before). By Lemma $3.4$ in \cite{GK}, for all $0\leq l\leq \nu $ ($\leq p+b-c$) we have
\begin{align*}
a_{j,l} &\equiv\begin{cases}
{{p+b-c-l}\choose b-m-j}{{c-1}\choose j}\quad \text{if}\ 0\leq j\leq b-m\\
{{p+b-c-l}\choose p+b-m-j}{{c-1}\choose j-1}\quad\text{if}\ b-m+1\leq j\leq c-m-1
\end{cases} \\
&\equiv \begin{cases}
{{p+b-c-l}\choose b-m-j}{{c-1}\choose j}\quad& \text{if}\ 0\leq j\leq b-m\\
0\quad&\text{if}\ b-m+1\leq j\leq c-m-1.
\end{cases}
\end{align*}
For each $m\in[b-\nu,\nu]\cap[1,b-1]$, a similar calculation to Case $(2)$ shows that $d_l\not\equiv0\bmod p$ for $l=b-m$ if $b\geq 2m+1$ and for $l=b-m-1$ if $b\leq 2m$.\\
\textbf{Subcase (II)} $b\leq m\leq \nu-1$\\
In this case, we observe that $\sigma_1(l) = 1$  for all $m\leq l\leq \nu-1$ as $\nu-1\leq p+b-c-1$. We claim that there exists at least one $l\in[m,\nu-1]$ such that $d_l\not\equiv 0(\bmod p)$. Note that $\epsilon_1 = 1$ as $b\leq 2m$, so $1\leq j\leq c-m-1$. By Lemma $3.5$ in \cite{GK}, we have
$$a_{j,l} \equiv (-1)^{l-m}\frac{{{p+b-m-1}\choose j-1}{{p-1+m-l}\choose c-1-j}}{{{p+b-m-c}\choose l-m}{{p+b-m-1}\choose c-1}}.$$
Using \eqref{case1_bino_rearr}, we have 
\begin{align*}
							d_l\equiv\frac{(-1)^l{{p-1+m-l}\choose m}}{{{p+b-m-c}\choose l-m}{{p+b-m-1}\choose c-1}}\underset{1\leq j\leq c-m-1}\sum{{p+b-m-1}\choose j-1}{{p-1-l}\choose c-m-1-j}.
						\end{align*}
Replace $j-1$ by $j$ in the above summation and by using Vandermonde's identity, we get
\begin{align*}
d_l\equiv\frac{(-1)^l{{p-1+m-l}\choose m}{{2p+b-m-2-l}\choose c-m-2}}{{{p+b-m-c}\choose l-m}{{p+b-m-1}\choose c-1}}.
\end{align*}
Now, observe that
\begin{align*}
{{2p+b-2m-2}\choose c-m-2}&=\frac{\underset{0\leq k\leq c-m-3}\Pi (2p+b-2m-2-k)}{(c-m-2)!}\\
&\equiv \frac{(-1)^{c-m}\underset{0\leq k\leq c-m-3}\Pi(2m+2-b+k)}{(c-m-2)!}\not\equiv0\bmod p
\end{align*}
because $0<2m+2-b+k\leq 2m+2-b+c-m-3\leq c+m-b-1\leq p-1$ (the last inequality follows as $2c-2-p\leq b$ \& $m\leq c-2$). Hence, we get $d_m\not\equiv0\bmod p$. Thus, for each $m\in[b,\nu-1]$ there exists at least one $l$ such that $d_l\not\equiv0\bmod p$.\\
\textbf{Case (4)} $2\leq b\leq 2(c-1)-p-1$ and $1 \leq \nu\leq p+b-c$\\
Note that $\epsilon_1\in\{0,1\}$ as $b\geq 2(p+b-c)+1-(p-1)\geq 2m+1-(p-1)$.\\
\textbf{Subcase (I)} $m\in[b-\nu,\nu]\cap[1,b-1]$\\
Observe that $\sigma_1(l)=0$ for all $l\in[0,\nu]$ (again taking $j = b-m$). We claim that there exists at least one $l\in[0,\nu]$ such that $d_l\not\equiv0\bmod\ p$. For each $m\in[b-\nu,\nu]\cap[1,b-1]$, using Lemma $3.4$ in \cite{GK}, we arrive at the same formula for $d_l$ as in Case (2). Thus, $d_l\not\equiv0\bmod p$ for $l=b-m$ if $b\geq 2m+1$ and for $l=b-m-1$ if $b\leq 2m$.\\
\textbf{Subcase (II)} $b\leq m\leq \nu-1$\\
We observe that $\sigma_1(l) = 1$ for all $m\leq l\leq \nu-1$ as $\nu-1\leq p+b-c-1$. For each $m\in[b,\nu-1]$, using Lemma $3.4$ in \cite{GK}, we have the same formula for $d_l$ as in Subcase (II) of Case (3), thus giving us that $d_l\not\equiv0\bmod p$ for $l=m$.\\
\textbf{Case (5)} $2\leq b\leq 2(c-1)-p-1$ and $p+b-c+1\leq\nu\leq c-3$\\
By Lemma $3.4$ in \cite{GK}, we have
\begin{eqnarray}\label{case 4 eqn}
\nonumber &{{r-l}\choose b-m+j(p-1)}\hspace{25em}\\
\equiv &\begin{cases}
{{p+b-c-l}\choose b-m-j}{{c-1}\choose j}\quad &\text{if}\ 0\leq l\leq p+b-c, \ 0\leq j\leq b-m\\
{{p+b-c-l}\choose p+b-m-j}{{c-1}\choose j-1}\quad&\text{if}\ 0\leq l\leq p+b-c,\ b-m+1\leq j\leq c-m-1\\
{{2p+b-c-l}\choose b-m-j}{{c-2}\choose j} \quad&\text{if}\ p+b-c+1\leq l\leq \nu,\ 0\leq j\leq b-m\\
{{2p+b-c-l}\choose p+b-m-j}{{c-2}\choose j-1}\quad&\text{if}\ p+b-c+1\leq l\leq \nu,\ b-m+1\leq j\leq c-m-1.
\end{cases}\hspace{2em}
\end{eqnarray}
\textbf{Subcase (I)} $1\leq m\leq b-1$\\
In this subcase, we observe that $\sigma_1(l) = 0$ for all $l\in[0,p+b-c]$ (take $j = b-m$). We claim that there exists at least one $l\in[0,p+b-c]\subset[0,\nu]$ such that $d_l\not\equiv 0(\bmod \ p)$. From \eqref{case 4 eqn} note that $a_{j,l}\equiv 0\bmod p$ if $b-m+1\leq j\leq c-m-1$ and $0\leq l\leq p+b-c$. Therefore, for $l\leq p+b-c$ we have using \eqref{case2_bino_rearr}
\begin{align*}
d_l&\equiv (-1)^m\underset{\epsilon_1\leq j\leq b-m}\sum {{p+b-c-l}\choose b-m-j}{{c-1}\choose j}{{c-1-j}\choose m}\\
&\equiv(-1)^m{{c-1}\choose m}\underset{\epsilon_1\leq j\leq b-m}\sum{{p+b-c-l}\choose b-m-j}{{c-m-1}\choose j}.
\end{align*}
By using Vandermonde's identity we get 
\begin{align*}
d_l\equiv \begin{cases}
(-1)^m{{c-1}\choose m}{{p+b-m-1-l}\choose b-m}\quad\text{if}\ b\geq 2m+1\\
(-1)^m{{c-1}\choose m}\left({{p+b-m-1-l}\choose b-m}-{{p+b-c-l}\choose b-m}\right)\quad\text{if}\ b\leq 2m.
\end{cases}
\end{align*}
Hence, if $b\geq 2m+1$, then for $l = b-m\leq \nu$ (as $\nu>p+b-c$) we have $d_l\equiv {{p-1}\choose b-m}\not\equiv 0\bmod p$. If $b\leq 2m$ then for $l =0$ we have $d_l\not\equiv 0\bmod p$ as ${{p+b-m-1}\choose b-m}\equiv {{b-m-1}\choose b-m} \equiv 0\bmod p$ and ${{p+b-c}\choose b-m}\not\equiv 0\bmod p$ (since $b-c<0$ and $p+b-c\geq b-m$). Thus, in both cases there exists at least one $l$ such that $d_l\not\equiv0\bmod p$.\\
\textbf{Subcase (II)} $b\leq m\leq p+b-c$\\
In this case, we observe that $\sigma_1(l) = 1$ for all $m\leq l\leq p+b-c$ and $\sigma_1(l) = 0$ for $l=p+b-c+1$ (take $j=1$). We claim that there exists at least one $l\in[m,p+b-c+1]$ such that $d_l\not\equiv 0(\bmod p)$. For $m\leq l\leq p+b-c$, Lemma $3.5$ in \cite{GK} gives
\begin{align*}
a_{j,l}\equiv \frac{(-1)^{l-m}{{p+b-m-1}\choose j-1}{{p-1+m-l}\choose c-1-j}}{{{p+b-m-c}\choose l-m}{{p+b-m-1}\choose c-1}}\bmod p.
\end{align*}
Using above and \eqref{case 4 eqn} for $a_{j,l}$ at $l=p+b-c+1$, we get 
\begin{align*}
d_l\equiv\begin{cases}
\beta(l)\underset{\epsilon_1\leq j\leq c-m-1}\sum {{p+b-m-1}\choose j-1}{{p-1-l}\choose c-m-1-j} \quad&\text{if}\ m\leq l\leq p+b-c\\
\beta(l)\underset{\epsilon_1\leq j\leq c-m-1}\sum {{2p+b-c-l}\choose p+b-m-j}{{c-m-2}\choose j-1}\quad&\text{if}\ l=p+b-c+1
\end{cases}
\end{align*}
where (using \eqref{case1_bino_rearr} and a variation of \eqref{case2_bino_rearr})
\[\beta(l) =\begin{cases}
\frac{(-1)^l{{p-1+m-l}\choose m}}{{{p+b-m-1}\choose c-1}{{p+b-m-c}\choose l-m}}\quad&\text{if}\ m\leq l\leq p+b-c\\
(-1)^m{{c-2}\choose m}\quad&\text{if}\ l=p+b-c+1. 
\end{cases}
\]
By replacing $j-1$ by $j$ we get 
\begin{align*}
d_l\equiv\begin{cases}
\beta(l)\underset{\epsilon_1-1\leq j\leq c-m-2}\sum {{p+b-m-1}\choose j}{{p-1-l}\choose c-m-2-j} \quad&\text{if}\ m\leq l\leq p+b-c\\
\beta(l)\underset{\epsilon_1-1\leq j\leq c-m-2}\sum {{2p+b-c-l}\choose p+b-m-1-j}{{c-m-2}\choose j}\quad&\text{if}\ l=p+b-c+1.
\end{cases}
\end{align*}
By using Vandermonde's identity we get
\begin{align*}
d_l\equiv \begin{cases}
\beta(l){{2p+b-m-2-l}\choose c-m-2}\quad&\text{if}\ m\leq l\leq p+b-c\ \&\ b\geq 2m+1-(p-1)\\
\beta(l)\left({{2p+b-m-2-l}\choose c-m-2}-{{p-1-l}\choose c-m-2}\right)\quad&\text{if}\ m\leq l\leq p+b-c\ \&\ b\leq 2m-(p-1)\\
\beta(l){{2p+b-m-2-l}\choose p+b-m-1}\quad &\text{if}\ l=p+b-c+1\ \&\ b\geq 2m+1-(p-1)\\
\beta(l)\left({{2p+b-m-2-l}\choose p+b-m-1}-{{2p+b-c-l}\choose p+b-m-1}\right)\quad &\text{if}\ l=p+b-c+1\ \&\ b\leq 2m-(p-1).\\
\end{cases}
\end{align*}
Hence, we get
\begin{align*}
d_l&\equiv\begin{cases}
\beta(l){{p+c-m-2}\choose c-m-2}\quad&\text{if}\ l =p+b-c\ \&\ b\geq 2m+1-(p-1)\\
\beta(l)\left({{p+c-m-3}\choose p+b-m-1}-{{p-1}\choose p+b-m-1}\right)\quad&\text{if}\ l =p+b-c+1\ \&\ b\leq 2m-(p-1).
\end{cases}\\
&\equiv\begin{cases}
\beta(l){{c-m-2}\choose c-m-2}\quad&\text{if}\ l =p+b-c\ \&\ b\geq 2m+1-(p-1)\\
-\beta(l){{p-1}\choose p+b-m-1}\quad&\text{if}\ l =p+b-c+1\ \&\ b\leq 2m-(p-1).
\end{cases}
\end{align*}
as ${{p+c-m-3}\choose p+b-m-1}\equiv {{c-m-3}\choose p+b-m-1}=0$ (note that $m\leq c-3$). From the above calculation, we observe that $d_l\not\equiv0\bmod p$ for $l = p+b-c$ if $b\geq 2m+1-(p-1)$ and for $l =p+b-c+1$ if $b\leq 2m-(p-1)$. Note that if $d_l\not\equiv0\bmod p$ for $l=p+b-c+1\leq\nu$, then $l\in[0, \nu-\sigma_1(l)]$ because $\sigma_1(p+b-c+1)=0$. Hence, in this case also there exists at least one $l$ such that $d_l\not\equiv0\bmod p$.\\
\textbf{Subcase (III)} $p+b-\nu-1\leq m\leq \nu$\\
Using $p+b-c+1\leq m\leq \nu$, we observe that $\sigma_1(l) = 0$ for all $p+b-c+1\leq l\leq \nu$ (take $j = p+b-m$). We claim that there exists at least one $l\in[p+b-c+1, \nu]$ such that $d_l\not\equiv0\bmod p$. For $p+b-c+1\leq l\leq \nu$ and $1\leq \epsilon_{1} \leq j\leq c-m-1$, using \eqref{case 4 eqn} we get
$a_{j,l}\equiv {{2p+b-c-l}\choose p+b-m-j}{{c-2}\choose j-1}.$
Therefore, we get
\begin{align*}
d_l&\equiv (-1)^m\underset{\epsilon_1\leq j\leq c-m-1}\sum {{2p+b-c-l}\choose p+b-m-j}{{c-2}\choose j-1}{{c-1-j}\choose m}\\
&\equiv(-1)^m{{c-2}\choose m}\underset{\epsilon_1\leq j\leq c-m-1}\sum {{2p+b-c-l}\choose p+b-m-j}{{c-m-2}\choose j-1}
\end{align*}
as ${{c-2}\choose j-1}{{c-1-j}\choose m} = {{c-2}\choose m}{{c-m-2}\choose j-1}$. By replacing $j-1$ by $j$ we get 
\begin{align*}
d_l&\equiv (-1)^m{{c-2}\choose m}\underset{\epsilon_1-1\leq j\leq c-m-2}\sum {{2p+b-c-l}\choose p+b-m-1-j}{{c-m-2}\choose j}\\
&\equiv (-1)^m{{c-2}\choose m}\underset{\epsilon_1-1\leq j\leq p+b-m-1}\sum {{2p+b-c-l}\choose p+b-m-1-j}{{c-m-2}\choose j}.
\end{align*}
The latter congruence follows as ${{c-m-2}\choose j} = 0$ for all $c-m-1\leq j\leq p+b-m-1$. By Vandermonde's identity, we get 
\begin{align*}
d_l\equiv\begin{cases}
(-1)^m{{c-2}\choose m}{{2p+b-m-2-l}\choose p+b-m-1}\quad&\text{if}\ b\geq 2m+1-(p-1)\\
(-1)^m{{c-2}\choose m}\left({{2p+b-m-2-l}\choose p+b-m-1}-{{2p+b-c-l}\choose p+b-m-1}\right)\quad&\text{if}\ b\leq 2m-(p-1).
\end{cases}
\end{align*}
Therefore, we have 
\begin{align*}
d_l&\equiv\begin{cases}
(-1)^m{{c-2}\choose m}{{p-1}\choose p+b-m-1}\ \ \ \ \ \ \text{if}\ l=p+b-m-1\ \&\ b\geq 2m+1-(p-1)\\
(-1)^m{{c-2}\choose m}\left({{p+c-m-3}\choose p+b-m-1}-{{p-1}\choose p+b-m-1}\right)\ \text{if}\ l =p+b-c+1\ \&\ b\leq 2m-(p-1)
\end{cases}\\
&\equiv \begin{cases}
(-1)^m{{c-2}\choose m}{{p-1}\choose p+b-m-1} &\text{if}\ l=p+b-m-1\ \&\ b\geq 2m+1-(p-1)\\
(-1)^{m+1}{{c-2}\choose m}{{p-1}\choose p+b-m-1}&\text{if}\ l =p+b-c+1\ \&\ b\leq 2m-(p-1).
\end{cases}
\end{align*} 
The latter congruence follows as ${{p+c-m-3}\choose p+b-m-1}\equiv {{c-m-3}\choose p+b-m-1}\equiv0\bmod p$ (note that $m\leq c-3$ and $c-m-3<p+b-m-1$), and ${{p-1}\choose p+b-m-1}\not\equiv0\bmod p$ as $b<m$. Hence, we have $d_l\not\equiv0\bmod p$ for $l = p+b-c+1$ if $b\leq 2m-(p-1)$ and for $l =p+b-m-1$ if $b\geq 2m+1-(p-1)$. Note that $l= p+b-m-1\leq \nu$ if $p+b-1-\nu\leq m$. Hence, for all $m\in[p+b-1-\nu, \nu]$ there exists at least one $l$ such that $d_l\not\equiv0\bmod p$.
\end{proof}

The following identities will only be required in \S \ref{pf_more_cases_conj}, where we prove more cases of Conjecture \ref{mat_conj}. 
\begin{lemma}{\label{mat cong}}
Let $b,i\in\mathbb{N}$ and suppose $c,l\in\mathbb{Z}_{\geq0}$ such that $l\geq b-c+1$ and $i\geq b-c+2$. Also assume that $1\leq m\leq c-1$ and $b\geq m$. Then
\begin{eqnarray*}
&\underset{0\leq k\leq c-m-1}\sum(-1)^k(c-k){{c-m-1+i}\choose k}{{c-m-2+i-k}\choose i-(b-c+2)}{{l+c-m-1-k}\choose l}\\
& = \begin{cases}
C_1+C_2 & \text{if}\ i\leq l\\
C_1+C_2+(-1)^{b-m+l}(m+1-i){{i-1}\choose l} & \text{if}\ l\leq i-1.
\end{cases}
\end{eqnarray*}
where 
$$C_1 =  \underset{0\leq j\leq  m_1}\sum (-1)^jc{{c-m-1+i}\choose i-(b-c+2+j)}{{l-(b-c+2+j)}\choose l-(b-m+1+j)}$$
$$C_2=\underset{0\leq j\leq m_2}\sum (-1)^j(b-m+1+j){{c-m-1+i}\choose i-(b-c+2+j)}{{l-(b-c+2+j)}\choose l-(b-m+j)}$$
where $m_1=\text{min}\{i-(b-c+2), l-(b-m+1)\}$ and $m_2=\text{min}\{i-(b-c+2), l-(b-m)\}$. Further, $C_1=0$ if $l \leq b-m$ and $C_2=0$ if $l \leq b-m-1$.
\end{lemma}
\begin{proof}
We obtain the lemma from Lemma \ref{app_C1_C2}, where we put $n = b-c+1$.    
\end{proof}

\begin{rk}\label{rk of mat cong}
If $b\geq 2m+1$ and $l\leq m$, then the above lemma gives $C_1 = 0= C_2$.
\end{rk}

\begin{lemma}\label{alpha' cong II}
Let $b,i\in\mathbb{N}$ and suppose $c,l\in\mathbb{Z}_{\geq0}$ such that $l\geq p+b-c+1$ and $i\geq p+b-c+2$. Also, assume that $1\leq m\leq c-2$ and $p+b-m-1\geq 0$. Then
	\begin{eqnarray*}
		&\underset{0\leq k\leq c-m-2}\sum(-1)^k(c-1-k){{c-m-2+i}\choose k}{{i+c-m-3-k}\choose i-(p+b-c+2)}{{l+c-m-2-k}\choose l}\\
		& = \begin{cases}
			C_1+C_2 & \text{if}\ l\geq i\\
			C_1+C_2+(-1)^{p+b-m+1+l}(m+1-i){{i-1}\choose l} & \text{if}\ l\leq i-1.
		\end{cases}
	\end{eqnarray*}
	where 
 $$C_1 = \underset{0\leq j\leq m_1}\sum (-1)^j(p+b-m+j){{c-m-2+i}\choose i-(p+b-c+2+j)}{{l-(p+b-c+2+j)}\choose l-(p+b-m-1+j)},$$
 \begin{align*}
		C_2 = \underset{0\leq j\leq m_2}\sum (-1)^j(c-1){{c-m-2+i}\choose i-(p+b-c+2+j)}{{l-(p+b-c+2+j)}\choose l-(p+b-m+j)}
\end{align*}
and $m_1=\text{min}\{i-(p+b-c+2), l-(p+b-m-1)\}$ and $m_2=\text{min}\{i-(p+b-c+2), l-(p+b-m)\}$. Further, $C_1 = 0$ if $l\leq p+b-m-2$ and $C_2 = 0$ if $l\leq p+b-m-1$. 
\end{lemma}
\begin{proof}
The lemma follows from Lemma \ref{app_C1_C2}, where we first replace $c$ with $c-1$ and then put $n = p+b-c+1$ (in this order). Next, we observe that $C_1$ \& $C_2$ of Lemma \ref{app_C1_C2} are $C_2$ \& $C_1$ respectively of the lemma above.
\end{proof}

\begin{rk}\label{rk of alpha' cong II}
	If $b\geq 2(m+1)-p$ and $l\leq m$ then above lemma gives $C_1 = 0= C_2$.
\end{rk}

\section{Determination of monomials in $\kerp$} \label{mono_kerp_section_II}

\begin{lemma}{\label{Vrm cong}}
Let $r,b,c,m\in\mathbb{Z}_{\geq 0}$ such that $2\leq b\leq p$ and $r\equiv b\bmod (p-1)$. Assume that $r\geq b+c(p-1)+1$. If\ \ $0\leq m\leq \text{min}\{c-1-\epsilon_1,p-1\}$, then for all $1\leq a\leq c-m-\epsilon_1$,
$$q(c-m-a)\equiv\underset{1\leq i\leq m+1}\sum(-1)^{i+1}\beta(a,i)q(c-m-1+i)\bmod(V^{(m+1)}_r).$$
\end{lemma}
\begin{proof}
Let $P_j= x^{r-(b+1)-(c-j+1)(p-1)}y^{b-2m-1+(c-m-j)(p-1)}$ for all $1\leq j\leq c-m-\epsilon_1$, and recall $q(i)= x^{r-b+m-i(p-1)}y^{b-m+i(p-1)}$ for all $n_0\leq i\leq c$, where $n_0 = 0$ if $b\geq m$ and $1$ otherwise. Indeed $P_j\in V_{r-(m+1)(p+1)}$ as $b-2m-1+(c-m-j)(p-1)\geq b-(2m+1)+\epsilon_1(p-1)\geq 0$.  Observe that 
\begin{align}{\label{thetap}}
&	\theta^{m+1}P_j = \underset{0\leq i\leq m+1}\sum (-1)^i{{m+1}\choose i} q(c-m-j+i).\hspace{2.8em}
\end{align}
We will prove the result by induction on $a$. Observe that for $a = 1$, the result is true, and can be seen by putting $j = 1$ in (\ref{thetap}) as follows
$$\underset{0\leq i\leq m+1}\sum (-1)^{i}{{m+1}\choose i}q(c-m-1+i)\equiv  0\bmod(V^{(m+1)}_r)$$
$$\implies\quad q(c-m-1)\equiv \underset{1\leq i\leq m+1}\sum (-1)^{i+1}{{m+1}\choose i}q(c-m-1+i)\bmod(V^{(m+1)}_r).$$
By induction step, assume the result is true for $1\leq a\leq k-1$, and we prove it for $a = k\leq c-m-\epsilon_1$.  By putting $j = k$ in \eqref{thetap} we have
$$ \underset{0\leq i\leq m+1}\sum (-1)^{i}{{m+1}\choose i}q(c-m-k+i)\equiv 0\bmod(V^{(m+1)}_r)$$
\begin{align*}
&\implies q(c-m-k) \equiv\underset{1\leq i\leq k-1}\sum(-1)^{i+1}{{m+1}\choose i}q(c-m-k+i)\\
&\hspace{10em}+\underset{k\leq i\leq m+1}\sum(-1)^{i+1}{{m+1}\choose i}q(c-m-k+i).
\end{align*}
We note that if $k>m+1$, then the second sum above is zero and in the first sum $i$ runs from $1$ to $m+1$ as ${{m+1}\choose i} = 0$ for all $m+1<i\leq k-1$.  We rename $i$ by $l$ in the first sum and replace $i-(k-1)$ by $i$ in the second sum. Thus, we get
\begin{align*}
&	q(c-m-k) \equiv\underset{1\leq l\leq k-1}\sum(-1)^{l+1}{{m+1}\choose l}q(c-m-k+l)\\
&	\hspace{10em}+\underset{1\leq i\leq m+1-(k-1)}\sum(-1)^{i+k}{{m+1}\choose i+k-1}q(c-m-1+i).
\end{align*}
Observe that the second sum can be taken over $1\leq i\leq m+1$ as ${{m+1}\choose i+k-1} = 0$ for all $i>m+1-(k-1)$. By induction, we have 
\begin{align*}
q(c-m-k) &\equiv \underset{1\leq l\leq k-1}\sum(-1)^{l+1}{{m+1}\choose l}\left(\underset{1\leq i\leq m+1}\sum(-1)^{i+1}\beta(k-l, i)q(c-m-1+i)\right)\\
&\hspace*{8em}+\underset{1\leq i\leq m+1}\sum(-1)^{i+k}{{m+1}\choose i+k-1}q(c-m-1+i)
\end{align*}
Now, we interchange the sums in the first sum and then combine them with the last sum, we get
\begin{align*}
q(c-m-a)\equiv \underset{1\leq i\leq m+1}\sum (-1)^{i+1}\beta(k,i)q(c-m-1+i)
\end{align*}
where $\beta(k,i)$ is defined in Lemma \ref{dfn beta}.
\end{proof}

\begin{cor}\label{fm from q(j)}
Let $r,b,c,m\in\mathbb{Z}_{\geq 0}$ such that $2\leq b\leq p$ and $r\equiv b\bmod (p-1)$. Assume that $r\geq b+c(p-1)+1$. Suppose  $1\leq m\leq \text{min}\{c-1-\epsilon_1,p-1\}$, and  $q(j)\in V_r^{(m+1)}+\kerp$ for all $c-m\leq j\leq c-1$. Then for all $\epsilon_1\leq j\leq c-m-1$, $$q(j)\equiv (-1)^{m}{{c-1-j}\choose m}q(c)(\bmod (\kerp + V_{r}^{(m+1)})).$$
\end{cor}

\begin{prop}{\label{gen 1}}
Let $r= s+p^t(p-1)d$, with $p\nmid d$, $s =b+c(p-1)$ and suppose also that $2\leq b\leq p, \ 1\leq c\leq p-1$.  Fix non zero slope $\nu(a_p)$ such that $0\leq m\leq\nu(a_p)< p-1$,  $(m,\nu)\not= (\nu, \nu(a_p))$, and also $s>2\nu(a_p)$. Further we assume  $t>\nu(a_p)+c-1$ if $(b,  c,  m) \not = (p,  1,  0)$ and $t>\nu(a_p)+c$ if $(b,  c,  m) = (p,  1,  0)$. Then for all $g\in G$ and for $0\leq l\leq \nu$ with $(l,\nu)\not= (\nu, \nu(a_p))$, there exists $f^l\in\indkg\left(\symqp\right)$ such that 
\begin{equation}{\label{1st gen 1}}
(T-a_p)f^l\equiv \left[g, \underset{\substack{0<j<s-m\\j\equiv (s-m) \bmod (p-1)}}\sum{{r-l}\choose j}x^{r-j}y^j\right].
\end{equation}
Further assume $(b,c, m)\not=(p, 1,0)$, $t>\nu(a_p)+c$ and $0\leq l,m\leq \nu-1$. If  $m = 0$, then assume $l$ further satisfies $\nu\left({{r-l}\choose p-1-l}\right)\geq 1$. If $l =\nu -1$ or $m = \nu-1$, then assume $\nu<\nu(a_p)$.  For all $g\in G$  and the above values of $l,m$ there exists $f^l\in\indkg\left( \symqp\right)$ such that  
\begin{equation}\label{2nd gen 1}
(T-a_p)\left(\frac{f^l}{p}\right)\equiv \left[g, \underset{\substack{0<j<s-m\\j\equiv (s-m) \bmod (p-1)}}\sum\frac{{{r-l}\choose j}}{p}x^{r-j}y^j\right].
\end{equation}
\end{prop}
\begin{proof}
    We begin by observing that the coefficients $\frac{{{r-l}\choose j}}{p}$ in \eqref{2nd gen 1} are integral if $0\leq m\leq l\leq b-c$ or $b\leq m\leq l\leq b-c+p$, and not integral in general.
						Consider the following functions
						\begin{eqnarray*}
							f_{3,l} &=& \sum_{\lambda\in I^*_1}\left[ g^0_{2,p\lambda}, \frac{F_l(x,y)}{\lambda^{m-l}p^l(p-1)}\right]\\
							f_{2,l} &=& \left[ g^0_{2,0},{{r-l}\choose r-m}\frac{F_m(x,y)}{p^m}\right]\\
							f_{1,l} &=&\left[g^0_{1,0}, \frac{1}{a_p}\underset{\substack{s-m\leq j<r-m\\j\equiv (r-m)\bmod (p-1)}}\sum{{r-l}\choose j}x^{r-j}y^j\right]\\
                                f_0 &=& \begin{cases}
							[1,\ F_s(x,y)] & \text{if}\ r\equiv m\bmod (p-1)\\
							0   \   & \    \text{else}.
						\end{cases}
						\end{eqnarray*}
                           Now,
						\begin{align*}
							&	T^+\left(\left[ g^0_{2,p\lambda}, \frac{F_l(x,y)}{\lambda^{m-l}p^l(p-1)}\right]\right)\\
							& = \sum_{\mu\in I^*_1}\left[ g^0_{3,p\lambda+p^2\mu}, \sum_{0\leq j\leq s-l}\frac{p^{j-l}(-\mu)^{s-l-j}}{\lambda^{m-l}(p-1)}\left({{r-l}\choose j}- {{s-l}\choose j}\right)x^{r-j}y^j\right]\\
							&	\hspace{4em}+\sum_{\mu\in I_1}\left[ g^0_{3,p\lambda+p^2\mu}, \sum_{s-l+1\leq j\leq r-l}\frac{p^{j-l}(-\mu)^{r-l-j}}{\lambda^{m-l}(p-1)}{{r-l}\choose j}x^{r-j}y^j\right]\\
							&	\hspace{15em} -\left[ g^0_{3,p\lambda},\frac{p^{s-2l}}{\lambda^{m-l}(p-1)}x^{r-s+l}y^{s-l}\right].
						\end{align*}
						\vspace{2mm}\\
						Now we will estimate the valuation of the coefficients in the three sums (I),(II) and (III) above.  In (I),  for $j\geq 1$,  $\nu\left({{r-l}\choose j}-{{s-l}\choose j}\right)\geq t-\nu(j!) \implies$ $ j-l+t-\nu(j!)\geq t-\nu+1>c\geq1$.  For (III),  $s-2l\geq s-2\nu>0$. For (II) the same computation as in (III) shows that $j-l>0$.  Therefore we have $T^+\left(f_{3,l}\right)\equiv\ 0\bmod p$. Note that $l\leq \nu-1$ in \eqref{2nd gen 1};  therefore, the same calculation gives the valuation of the coefficients to be positive as $s-2l-1\geq s-2\nu+1>1$ and for $j\geq 1$, $j-l-1+t-\nu(j!)\geq t-(\nu-1)>c> 0$. Hence we have $T^+\left(\frac{f_{3,l}}{p}\right)\equiv\ 0\bmod p$. Now, 
						\begin{align*}
							T^-\left(\left[ g^0_{2,p\lambda}, \frac{F_l(x,y)}{\lambda^{m-l}p^l(p-1)}\right]\right) =- \left[ g^0_{1,0}, \sum_{0\leq j\leq s-l}\frac{p^{r-s}\lambda^{s-m-j}}{(p-1)}{{s-l}\choose j}x^{r-j}y^j\right]\\
							+\left[ g^0_{1,0}, \sum_{0\leq j\leq r-l}\frac{\lambda^{r-m-j}}{(p-1)}{{r-l}\choose j}x^{r-j}y^j\right].
						\end{align*}
						In the first sum, the valuation of the coefficients is at least $r-s\gg 0$, and therefore we have 
						\begin{eqnarray*}
							T^-\left(\frac{ f_{3,l}}{p}\right) &\equiv & \left[g^0_{1,0}, \underset{\substack{0\leq j\leq r-m\\ j\equiv (r-m)\bmod (p-1)}}\sum\frac{{{r-l}\choose j}}{p}x^{r-j}y^j\right].\\
						\end{eqnarray*}
						\begin{align*}
							 &T^+(f_{2,l}) = \sum_{\mu\in I^*_1}\left[ g^0_{3,p^2\mu}, \sum_{0\leq j\leq s-m} p^{j-m}(-\mu)^{s-m-j}{{r-l}\choose r-m}\left({{r-m}\choose j}-{{s-m}\choose j}\right) x^{r-j}y^j\right]\\
							& \hspace{9em}+\sum_{\mu\in I_1}\left[ g^0_{3,p^2\mu}, \sum_{s-m+1\leq j\leq r-m}\frac{p^j(-\mu)^{r-m-j}}{p^m}{{r-l}\choose r-m}{{r-m}\choose j}x^{r-j}y^j\right]\\
							&\hspace{22em}  -\left[g^0_{3,0},\ p^{s-2m}{{r-l}\choose r-m}x^{r-s+m}y^{s-m}\right].
						\end{align*}
						In (I) sum,  for $j\geq 1$,  $ j-m+t-\nu(j!)\geq t-\nu+1>c\geq1$.  For (III),  $s-2m\geq s-2\nu>0$. For (II), the same computation as in (III) shows that $j-m>0$. Therefore we have $T^+(f_{2,l})\equiv\ 0\bmod p$. Note that $m\leq \nu-1$ in \eqref{2nd gen 1};  therefore, the same calculation gives the valuation of the coefficients to be positive as $s-2m-1\geq s-2\nu+1>1$ and for $j\geq 1$, $j-m-1+t-\nu(j!)\geq t-(\nu-1)>c> 0$. Hence we have $T^+\left(\frac{f_{2,l}}{p}\right)\equiv\ 0\bmod p$. Now, 
						\begin{align*}
							T^-(f_{2,l}) &= -\left[g^0_{1,0}, p^{r-s}{{r-l}\choose r-m}x^{r-s+m}y^{s-m}\right]+\left[g^0_{1,0}, {{r-l}\choose r-m}x^my^{r-m}\right]\\
							\implies\ \ T^-(f_{2,l}) &\equiv \left[g^0_{1,0},{{r-l}\choose r-m}x^my^{r-m}\right] \quad (\text{as}\quad r-s\gg 0)\\
							\& \ \ \ \ T^-\left(\frac{ f_{2,l}}{p}\right)  &\equiv \left[g^0_{1,0},\frac{{{r-l}\choose r-m}}{p}x^m y^{r-m}\right] \quad (\text{as}\quad r-s-1\gg 0).
						\end{align*}
						Now,
                             \begin{align}\label{tf1}
                                 \nonumber & T^+(f_{1,l}) = \sum_{\lambda\in I^*_1}\left[ g^0_{2,p\lambda}, \sum_{0\leq j\leq r}\frac{p^j(-\lambda)^{s-m-j}}{a_p}\underset{\substack{s-m\leq i<r-m\\i\equiv (r-m)\bmod (p-1)}}\sum{{r-l}\choose i}{i\choose j}x^{r-j}y^j\right]\\
							&\hspace{16em} +\left[ g^0_{2,0}, \underset{\substack{s-m\leq j<r-m\\j\equiv (r-m)\bmod (p-1)}}\sum\frac{p^j}{a_p}{{r-l}\choose j}x^{r-j}y^j\right]. 
                             \end{align}
						\vspace{2mm}\\
					 Here we note that $s-(\nu+1)\leq \text{min}\{s-l-1,s-m\}$ as $l,m\leq \nu$. Observe that the first sum truncates to $j\leq s-(\nu+1)$ and the second sum is zero mod $p$ as $j-\nu(a_p)>0$ for all $j\geq s-\nu$. Hence
						\begin{align*}
							T^+(f_{1,l}) = \sum_{\lambda\in I^*_1}\left[ g^0_{2,p\lambda}, \sum_{0\leq j\leq s-(\nu+1)}\frac{p^j(-\lambda)^{s-m-j}}{a_p}\underset{\substack{s-m\leq i<r-m\\i\equiv (r-m)\bmod (p-1)}} \sum{{r-l}\choose i}{i\choose j}x^{r-j}y^j\right]
						\end{align*}
						$$\implies T^+(f_{1,l}) = \sum_{\lambda\in I^*_1}\left[ g^0_{2,p\lambda}, \sum_{0\leq j\leq s-(\nu+1)}\frac{p^j(-\lambda)^{s-m-j}}{a_p}S_{r,j,l,m}\ x^{r-j}y^j\right]$$
      where $S_{r,j,l,m}$ is defined in equation $(3.1)$ in \cite{GK}.  Now $(b,c,m)\not=(p,1,0)$ implies that either $b\leq p-1$ or $c+m\geq 2$ (or both).  So Lemma $3.3 (II)$ in \cite{GK} gives $\nu(S_{r,j,l,m})\geq t+1-c$,  and therefore the valuation of above coefficients is at least $j+t+1-c-\nu(a_p)\geq t-(\nu(a_p)+c-1)>0$ giving $ T^+(f_{1,l})\equiv 0 \bmod p$. For $(b,c,m)= (p,1,0)$, Lemma $3.3 (II)$ in \cite{GK} gives $\nu(S_{r,j,l,m})\geq t-c$, so the valuation of above coefficients is at least $j+t-c-\nu(a_p)\geq t-(c+\nu(a_p))>0$, hence $T^+(f_{1,l})\equiv\ 0\bmod p$.\\
					   Again note that $s-\nu\leq \text{min}\{s-l-1,s-m\}$ as $l,m\leq \nu-1$ in \eqref{2nd gen 1}. Note also that the first sum truncates to $j\leq s-\nu$, and the second sum is zero mod $p$ as $j-\nu(a_p)-1>0$ for all $j\geq s-\nu+1$. Hence,
							$$T^+\left(\frac{f_{1,l}}{p}\right) = \sum_{\lambda\in I^*_1}\left[ g^0_{2,p\lambda}, \sum_{0\leq j\leq s-\nu}\frac{p^{j-1}(-\lambda)^{s-m-j}}{a_p}S_{r,j,l,m}\ x^{r-j}y^j\right]$$
							Since in this case $(b,c,m)\not=(p,1,0)$, Lemma $3.3 (II)$ in \cite{GK} gives $\nu(S_{r,j,l,m})\geq t+1-c$ for all $j\leq s-\nu$,  and therefore the valuation of above coefficients is at least $j-1+t+1-c-\nu(a_p)\geq t-(\nu(a_p)+c)>0$ giving $ T^+(\frac{f_{1,l}}{p})\equiv\ 0\bmod p$. 
						Next,
						$$T^-(f_{1,l}) = \left[ 1, \underset{\substack{s-m\leq j<r-m\\ j\equiv (r-m)\bmod (p-1)}}\sum\frac{p^{r-j}}{a_p}{{r-l}\choose j}x^{r-j}y^j\right]$$
						where the valuation of the coefficients is at least $r-j-\nu(a_p)\geq m+p-1-\nu(a_p)>0\ \ \Rightarrow\ \ T^-(f_{1,l})\equiv\ 0\bmod\ p$.  Also, 
						$$T^-\left(\frac{f_{1,l}}{p}\right) = \left[ 1, \underset{\substack{s-m\leq j<r-m\\ j\equiv (r-m)\bmod (p-1)}}\sum\frac{p^{r-j-1}}{a_p}{{r-l}\choose j}x^{r-j}y^j\right].$$
						For $j\leq r-m-2(p-1)$ the valuation of the above coefficients are at least $r-j-1-\nu(a_p)\geq p-2+m+p-1-\nu(a_p)>0$.  For $j= r-m-(p-1)$, the valuation of the coefficient is
						$r-j-1-\nu(a_p)+\nu({{r-l}\choose r-m-(p-1)})\geq m+\nu\left({{r-l}\choose r-m-(p-1)}\right)-1+p-1-\nu(a_p)>m+\nu\left({{r-l}\choose r-m-(p-1)}\right)-1\geq 0$. Observe that the last inequality is clear as either $m\geq 1$ or $\nu{{r-l}\choose p-1-l}\geq 1$. 
						Therefore we have $T^-\left(\frac{f_{1,l}}{p}\right) \equiv 0\bmod p$.
      
					      If $r\equiv m\bmod (p-1)$, then
							\begin{eqnarray*}
								T^+(f_0) &=& \sum_{\lambda\in I^*_1}\left[ g^0_{1,\lambda}, (-1+(-\lambda)^{r-s})x^r\right]+[g^0_{1,0},-x^r]\\
								&&	+\sum_{\lambda\in I_1}\left[g^0_{1,\lambda},\sum_{1\leq j\leq r-s}p^j{{r-s}\choose j}(-\lambda)^{r-s-j}x^{r-j}y^j\right]\\
								\Rightarrow \ \ \ T^+(f_0) &\equiv & -[g^0_{1,0},\  x^r]\\
								T^-(f_0) &=& \left[ \alpha, \ -p^rx^r + p^sx^sy^{r-s}\right]\equiv\ 0\bmod p\\
								T^+\left(\frac{f_0}{p}\right) &=& \sum_{\lambda\in I^*_1}\left[ g^0_{1,\lambda}, \frac{(-1+(-\lambda)^{r-s})}{p}x^r\right]+[g^0_{1,0},\ -\frac{1}{p}x^r]\\
								&&+\sum_{\lambda\in I_1}\left[g^0_{1,\lambda},\sum_{1\leq j\leq r-s}p^{j-1}{{r-s}\choose j}(-\lambda)^{r-s-j}x^{r-j}y^j\right]
							\end{eqnarray*}
							Observe that if $j = 1 $, $\ {{r-s}\choose j} = r-s$ which is divisible  by $p^t,\ t\geq 1$. Thus
							\begin{eqnarray*}
								\ \ \ T^+\left(\frac{f_0}{p}\right) &\equiv & -[g^0_{1,0},\ \frac{1}{p} x^r]\\
								T^-\left(\frac{f_0}{p}\right) &=& \left[ \alpha, \ -p^{r-1}x^r + p^{s-1}x^sy^{r-s}\right]\equiv\ 0\bmod p
						\end{eqnarray*} 
						 Note that $a_pf_{3,l},a_pf_{2,l},a_pf_{0}\equiv  0$, as $l,m<\nu(a_p)$ from the hypotheses of \eqref{1st gen 1}. Similarly, $\frac{a_pf_{3,l}}{p},\frac{a_pf_{2,l}}{p}$ since $l,m<\nu(a_p)-1$ from the hypotheses of \eqref{2nd gen 1}. Note that $\frac{a_pf_{0}}{p}\equiv 0$ since $\nu(a_p)>1$. It is clear from the hypotheses that $\nu(a_p) = 1$ is not allowed for \eqref{2nd gen 1}.
						\begin{eqnarray*}
							(T-a_p)(f_{3,l})&\equiv &\left[g^0_{1,0}, \underset{\substack{0\leq j\leq r-m\\j\equiv (r-m)\bmod (p-1)}}\sum{{r-l}\choose j}x^{r-j}y^j\right]\\
							(T-a_p)(f_{2,l})&\equiv &\left[g^0_{1,0},{{r-l}\choose r-m}x^my^{r-m}\right]\\
							(T-a_p)(f_{1,l})&\equiv &-\left[g^0_{1,0}, \underset{\substack{s-m\leq j<r-m\\j\equiv (r-m)\bmod (p-1)}}\sum{{r-l}\choose j}x^{r-j}y^j\right]
						\end{eqnarray*}
						\[(T-a_p)(f_{0,l})\equiv  \begin{cases}
							-[g^0_{1,0},\  x^r] & \ if\ r\equiv m\bmod (p-1)\\
							0   \   & \    else
						\end{cases}
						\]
						Hence  $f^l := f_{3,l}-f_{2,l}+f_{1,l}+f_{0,l}$ gives the required result. Finally, we deduce \eqref{1st gen 1} and \eqref{2nd gen 1} from the fact that $T$ is $G$-linear, and so we can replace $g_{1,0}^0$ by any $g\in G$.
\end{proof}

\begin{prop}{\label{gen m 1}}
Let $r= s+p^t(p-1)d$, with $p\nmid d$, $s =b+c(p-1)$ and also suppose that $2\leq b\leq p, \ 1\leq c\leq p-1$.  Fix $a_p$ such that $1<\nu(a_p)<p-1$ and $s>2\nu(a_p)$.   Assume $1\leq m\leq\text{min}\{\nu, c-1-\epsilon_1\}$ and $(m,\nu)\not= (\nu, \nu(a_p))$. If  $t>\nu(a_p)+c-1$, then for all $g\in G$ and for $0\leq l\leq \nu$ with $(l,\nu)\not= (\nu, \nu(a_p))$, there exists $f^l\in \indkg\left(\symqp\right)$ such that 
\begin{equation}{\label{1st gen m 1}}
(T-a_p)f^l\equiv \left[g, \underset{1\leq i\leq m+1}\sum\alpha(i,l)q(i+c-m-1)\right](\bmod X)
\end{equation}
where $\alpha(i,l)$ is as in \eqref{dfn alpha} and $X:= V^{(m+1)}_r+\kerp$. Further, suppose $1\leq m\leq\text{min}\{\nu-1, c-1-\epsilon_1\}$, $0\leq l\leq \nu-1$, and $t>\nu(a_p)+c$. Assume $\nu<\nu(a_p)$ if $l =\nu -1$ or $m = \nu-1$. If for some $l$ as above ${{r-l}\choose b-m+j(p-1)}\equiv 0\bmod p$ for $0\leq j\leq c-1$ such that $b-m+j(p-1)>0$, then for all $g\in G$ there exists $f^l\in \indkg\left(\symqp\right)$ such that  
\begin{equation}{\label{2nd gen m 1}}
(T-a_p)\left(\frac{f^l}{p}\right)\equiv \left[g, \underset{1\leq i\leq m+1}\sum\frac{\alpha(i,l)}{p}q(i+c-m-1)\right](\bmod X).
\end{equation}
\end{prop}
\begin{proof}
First, we note that 
\begin{align*}
&	A:= \underset{\substack{0<j<s-m\\j\equiv (s-m) \bmod (p-1)}}\sum{{r-l}\choose j}x^{r-j}y^j\\
&	 = \underset{\substack{0\leq j'\leq c-1\\ b-m+j'(p-1)>0}}\sum{{r-l}\choose b-m+j'(p-1)}q(j').
\end{align*}
Now, we write the last sum into three parts: $0\leq j<\epsilon_1$,  $\epsilon_1\leq  j'\leq c-m-1$, and $c-m\leq j'\leq c-1$.  Here we note that Remark $4.4$ in \cite{KBG} gives the first sum belongs to $\kerp$ as $m< \nu(a_p)$.  In the second range of sum, we put $a = c-m-j'$ and in the third range of sum we put $i' = j'-(c-m-1)$. Therefore, we have 
\begin{align*}
&A \equiv  \underset{1\leq a\leq c-m-\epsilon_1}\sum{{r-l}\choose b-m+(c-m-a)(p-1)}q(c-m-a)\\
&\hspace{2em}+ \underset{1\leq i'\leq m}\sum{{r-l}\choose b-m+(c-m-1+i')(p-1)}q(c-m-1+i')\bmod (\kerp).
\end{align*}
By Lemma \ref{Vrm cong}, for $1\leq a\leq c-m-\epsilon_1$ we have
$$q(c-m-a)\equiv\underset{1\leq i\leq m+1}\sum(-1)^{i+1}\beta(a,i)q(c-m-1+i)\bmod(V^{(m+1)}_r).$$
Therefore, we get 
\begin{align*}
&A\equiv  \underset{1\leq i'\leq m}\sum{{r-l}\choose b-m+(c-m-1+i')(p-1)}q(c-m-1+i')\\ &+\underset{1\leq a\leq c-m-\epsilon_1}\sum{{r-l}\choose b-m+(c-m-a)(p-1)}\left(\underset{1\leq i\leq m+1}\sum(-1)^{i+1}\beta(a,i)q(c-m-1+i)\right).
\end{align*}
The above congruence is over $\bmod (V_r^{(m+1)}+\kerp)$. Now, we interchange the sums in the last sum and combine them with the first sum (replace $i'$ by $i$); we get
\begin{align*}
&	A \equiv \underset{1\leq i\leq m+1}\sum \alpha(i,l)q(c-m-1+i)\bmod (V^{m+1}_r+\kerp)
\end{align*}
where $\alpha(i,l)$ is as in \eqref{dfn alpha}. Hence, we obtained the first part of our result just by using Proposition \ref{1st gen 1}. For the second part, we rerun the proof with
$A:= \underset{\substack{0<j<s-m\\j\equiv (s-m) \bmod (p-1)}}\sum\frac{{{r-l}\choose j}}{p}x^{r-j}y^j$
by noting that $\frac{{{r-l}\choose j}}{p}$ is integral.
\end{proof}

\begin{prop}\label{Fm_kerp_conj}
Let $r= s+p^t(p-1)d$, with $p\nmid d$, $s =b+c(p-1)$ and also suppose that $2\leq b\leq p, \ 1\leq c\leq p-1$.  Fix $a_p$ such that $\nu(a_p)$ is non-integral, $s>2\nu(a_p)$, and let $1\leq m\leq\nu\leq c-1-\epsilon$. Also assume  $t>\nu(a_p)+c$, and  $q(j)\in (V_r^{(m+1)}+\kerp$ for all $c-m\leq j\leq c-1$. Then $F_m(x,y)\in(V_r^{(m+1)}+\kerp)$ in each of the following cases (with the additional condition on $m$)given below:
\begin{enumerate}
\item $1\leq m\leq \text{min}\{\nu-1,b-c\}$ if $c-1\leq b\leq p$.
\item  $b-\nu\leq m\leq \nu$ if $c-1\leq b\leq p$.
\item $m\in\left([b-\nu,\nu]\cap[1,b-1]\right)\cup[b,\nu-1]$ if $2\leq b\leq c-2$ and $1\leq \nu\leq\text{min}\{c-2,p+b-c\}$.
\item $m\in[1,p+b-c]\cup[p+b-\nu-1, \nu]$ if $2\leq b\leq 2(c-1)-p-2$ and $p+b-c+1\leq \nu\leq c-3$.
\end{enumerate}
\end{prop}
\begin{rk}
(i) We note that the hypothesis $q(j)\in (V_r^{(m+1)}+\kerp$ for all $c-m\leq j\leq c-1$ (see Conjecture \ref{mono II conj}) in Proposition \ref{Fm_kerp_conj} is crucial for applying Corollary \ref{fm from q(j)}.\\
(ii) Note that the statement $(3)$ of the proposition above is obtained by combining $(3)$ and $(4)$ of Lemma \ref{d_l_non_zero}.\\
(iii) Corollary \ref{fm from q(j)} is applicable since $c-1-\epsilon \leq c-1-\epsilon_{1}$ if $1\leq m \leq c-1-\epsilon$.
\end{rk}
\begin{proof}
Recall the definition of $\sigma_1(l)$ from \eqref{sigma_1l} and observe that $\frac{{{r-l}\choose b-m+j(p-1)}}{p^{\sigma_1(l)}}$ is\vspace{.5em}\\ integral. By Remark 4.4 in \cite{KBG}, we have $q(j)\in \kerp$ for all $j\in[0,\epsilon_1-1]$ such that $b-m+j(p-1)>0$. Therefore, for each $l\in[0, \nu-\sigma_1(l)]$, Proposition \ref{gen 1} gives
$$(T-a_p)\left(\frac{f^l}{p^{\sigma_1(l)}}\right)\equiv \left[g,\ \underset{\epsilon_1\leq j\leq c-m-1}\sum \frac{{{r-l}\choose b-m+j(p-1)}}{p^{\sigma_1(l)}}q(j)\right]\bmod (V_r^{(m+1)}+\kerp)$$
as $q(j)\in (V_r^{(m+1)}+\kerp$ for all $j\in[c-m, c-1]$ (by the hypothesis). Since $\nu (a_p)$ is non-integral, we can also apply Proposition \ref{gen 1} in the case when $l = \nu$ in \eqref{1st gen 1} or when $l, m = \nu-1$ in \eqref{2nd gen 1}.
By Corollary \ref{fm from q(j)}, we have
\begin{align*}
(T-a_p)\left(\frac{f^l}{p^{\sigma_1(l)}}\right)&\equiv\left[g,\ d_l\  q(c)\right]\bmod (\kerp+V_r^{(m+1)}),
\end{align*}
where $d_l$ is as in \eqref{dfn_d_l}. For each $b,c,\nu$ and $m$ as in the statement, Lemma \ref{d_l_non_zero} gives that there exists at least one $l\in[0,\nu-\sigma_1(l)]$ such that $d_l\not=0\bmod p$. Hence, we get $q(c)\in(\kerp + V_r^{(m+1)})$. Further, note that Remark 4.4 in \cite{KBG} gives $x^my^{r-m}\in\kerp$ as $m<\nu(a_p)$. Consequently, $F_m(x,y)\equiv q(c)\bmod(\kerp + V_r^{(m+1)})$, thereby giving our proposition.
\end{proof}

\section{Conjecture}\label{conje_local_sec}
We consider the following matrix
\begin{equation}\label{matrix_A_conj}
A = \left(\alpha(i,l)\right)_{\substack{1\leq i\leq m+1\\ 0\leq l\leq m}}, 
\end{equation} 
where $\alpha(i,l)$ are defined in \eqref{dfn alpha}. We note that $\alpha(i,l)$ has a complicated formula, so proving that the linear system of equations $AX=e_i$ has a solution is generally hard. However, we make the following conjecture about the solutions to the above linear systems based on SageMath programs. 
					
\begin{conj}[Matrix form]\label{mat_conj}
Let $r= s+p^t(p-1)d$, with $p\nmid d$, $s =b+c(p-1)$ and suppose that $2\leq b\leq p, \ 1\leq c\leq p-1$, and $t\geq 2$. Suppose also that $1\leq m\leq c-1-\epsilon$. Then the linear systems $AX=e_i(\bmod p)$ has a solution for all $1\leq i\leq m$.
\end{conj}

We note that if $AX = e_{i'}(\bmod p)$ has a solution, then under the hypotheses of Proposition \ref{gen m 1}, we get $q(c-m-1+i')\in(V_r^{(m+1)}+\kerp)$. One can see this as follows, if $X^t = (d_{0},d_{1},...,d_{m})$ be a solution of $AX=e_{i'}(\bmod p)$, then take $f:=\underset{0\leq l\leq m}\sum d_lf^l$, where $f^l$ are defined in Proposition \ref{gen m 1}. Observe that
\begin{align*}
(T-a_p)f& \equiv \left[g,\underset{1\leq i\leq m+1}\sum\left(\underset{0\leq l\leq m}\sum d_l\alpha(i,l)\right)q(c-m-1+i)\right]\bmod (V_r^{(m+1)}+\kerp)\\
& \equiv \left[g,\ q(c-m-1+i')\right]\bmod (V_r^{(m+1)}+\kerp).
\end{align*}
Therefore, as a consequence of the above conjecture, we get the following conjecture for monomials $q(j)$.
\begin{conj}[Monomial form]\label{mono II conj}
Let $r= s+p^t(p-1)d$, with $p\nmid d$, $s =b+c(p-1)$ and suppose that $2\leq b\leq p, \ 1\leq c\leq p-1$.  Fix $a_p$ such that $1<\nu(a_p)<p-1$ and $t>\nu(a_p)+c$. Suppose also that $s>2\nu(a_p)$. Further we assume that $1\leq m\leq \text{min}\{\nu,c-1-\epsilon\}$ and $(m,\nu)\not=(\nu,\nu(a_p))$. Then the monomials $q(j)\in(V_r^{(m+1)}+\kerp)$ for all $c-m\leq j\leq c-1$.
\end{conj}

\begin{rk}
    Note that the crucial hypothesis in Proposition \ref{Fm_kerp_conj} on $q(j)$ is precisely the claim of Conjecture \ref{mono II conj}. 
\end{rk}

In the following lemma, we prove enough cases of Conjecture \ref{mat_conj} so that Proposition \ref{Fm_kerp_conj} $(1)$ holds unconditionally. In the last section, we discuss evidence for the remaining cases (and also provide proof in many cases) of Conjecture \ref{mat_conj}.
\begin{lemma}\label{proof_conj_case1}
If $c-1\leq b\leq p$ and $1\leq m\leq \min\{b-c,c-1-\epsilon\}$, then Conjecture \ref{mat_conj} is true.
\end{lemma}

\begin{proof}
First, we note that the range of $m$ is non-empty only when $b\geq c+1$, so we will prove our lemma only for $b\geq c+1$.
Now, let's express $A$ in \eqref{matrix_A_conj} as follows 
\begin{align}\label{pf_conj_case_1_eqn}
A = \begin{pmatrix}
     A'&B'\\
     A''&B''
     \end{pmatrix},
\end{align}
where the ranges of $i$ and $l$ are divided into non-empty intervals $[1,m], [m+1,m+1]$ and $[0,m-1]$, $[m,m]$, determining the order of the blocks. For $0\leq l\leq m$ and $1\leq i\leq m+1$, Lemma \ref{alpha cong} gives $\alpha_1(i,l)\equiv0\bmod p$ (as $m\leq b-c$), and also that the additional term of $\alpha (i,l)$ (for $i \leq m$) is given by 
$$X(i,l) = { {r-l}\choose b-m+(i+c-m-1)(p-1)}\equiv {{b-c-l}\choose b-c+1-i}{c\choose i+c-m-1}.$$ Hence, we get
$$\alpha(i,l)\equiv \begin{cases}{{b-c-l}\choose b-c+1-i}{c\choose c-m-1+i} & \text{if}\ 0\leq l\leq m, 1\leq i\leq m\\
 0 & \text{if}\ 0\leq l\leq m, i = m+1.
\end{cases}$$
The above congruence implies that modulo $p$, $A'$ is an invertible lower triangular matrix (with all the diagonal entries non-zero modulo $p$, given by $i = l+1$), and $A'', B''$ are zero modulo $p$. Hence, for every $1\leq i\leq m$, modulo $p$ the row rank of $[A:e_{i'}]$ is same as the row rank of $A$. Thus, the linear systems $AX\equiv e_{i'}(\bmod \ p)$ has a solution for all $1\leq i'\leq m$.

\end{proof}

 \section{Elimination of JH factor}

\begin{prop}\label{rm_JH_conj}
Let $r= s+p^t(p-1)d$, with $p\nmid d$, $s =b+c(p-1)$ and suppose also that $2\leq b\leq p, \ 2\leq c\leq p-1$.  Fix $a_p$ such that $\nu(a_p)$ is non-integral, $1<\nu(a_p)<c-\epsilon$ and $t> \nu(a_p)+c$. Further, assume \textbf{Conjecture \ref{mono II conj}} is true. Then the map $P$ surjects from $\indkg\left(\frac{V_r^{(n_1+1)}}{V_r^{(n_2)}}\right)$, where $n_1$ and $n_2$ are defined as follows:
\begin{enumerate}
    \item If $c-1\leq b\leq p$, then $n_1 = \text{min}\{b-c,\nu-1\}$ and $n_2 = \text{min}\{b-\nu, \nu+1\}$.
    \item Suppose $2\leq b\leq c-2$ and $1\leq \nu\leq \text{min}\{c-2,p+b-c\}$. Then
    $$n_1 = \begin{cases}
        \nu-1 \quad \text{if}\ b\leq \nu\\
        -1 \quad \text{otherwise}
    \end{cases}
    \quad \&\quad n_2 = \begin{cases}
        \nu+1 \quad\text{if} \ b\leq \nu\\
        \text{min}\{b-\nu,\ \nu+1\}\quad\text{otherwise}.
    \end{cases}$$
    \item If $2\leq b\leq 2c-4-p$ and $p+b-c+1\leq\nu\leq c-3$, then $n_1 = p+b-c$ and  $n_2 =\text{min}\{p+b-\nu-1,\ \nu+1\}$.
\end{enumerate}
\end{prop}
\begin{rk}\label{conj_required_1}
   In the proposition above, Conjecture \ref{mono II conj} is needed to apply Proposition \ref{Fm_kerp_conj}. Note that the statement $(1)$ above is obtained by combining the statements $(1)$ and $(2)$ of Proposition \ref{Fm_kerp_conj}. If $b\geq c+\nu-1$, then the range of $m$ in Proposition \ref{Fm_kerp_conj} $(2)$ is empty. Therefore, the case $b\geq c+\nu-1$ (appearing later) is covered entirely by the statement of Proposition \ref{Fm_kerp_conj} $(1)$.  
\end{rk}
\begin{proof}
In all three parts of the proposition, we need to show that the map $P$ factors through  $\indkg\left(\frac{V_r^{(n_1+1)}}{V_r^{(n_2)}}\right)$. In order to do that we apply Lemma \ref{lm m<c} \& \ref{lm m>c} to the intervals $[0,n_1]$ \& $[n_2,\nu]$ respectively, whenever they are non-empty. To apply Lemma \ref{lm m>c} in each cases, we show that $F_m(x,y)\in\kerp$ for all $m\in [n_2,\nu]$. In fact, it suffices to show that $F_m(x,y)\in(V_r^{(m+1)}+\kerp)$ for all $m\in[n_2,\nu]$. We show this by reverse induction on $m$ after observing that $F_{\nu}(x,y)\in\kerp$ because $F_{\nu}(x,y)\in(V_r^{(m+1)}+\kerp)$ and $\indkg(V_r^{(\nu+1)})\subset\kerp$.

Next, we explain the steps to apply Lemma \ref{lm m<c} for $[0,n_1]$. If $F_m(x,y)\in(V_r^{(m+1)}+\kerp)$ for $m\geq 1$, then there exists $v_{m+1}\in V_r^{(m+1)}$ such that $F_m(x,y)-v_{m+1}\in\kerp$. Let $W_m$ be the submodule of $V_r^{(m)}$ generated by $F_m(x,y)-v_{m+1}$. By using Lemma \ref{lm3.2}, and given that $F_m(x,y)-v_{m+1}\in\kerp$, we observe that $W_m$ satisfies the required conditions of Lemma \ref{lm m<c}. Hence, to apply Lemma \ref{lm m<c} for $[0,n_1]$, we show that $F_m(x,y)\in(V_r^{(m+1)}+\kerp)$ for all $m\in[1,n_1]$ and construct the required $W_0$ separately. In each of the above cases, Proposition \ref{Fm_kerp_conj} gives $F_m(x,y)\in(V_r^{(m+1)}+\kerp)$ for all $m\in[1,n_1]\cup [n_2,\nu]$. Therefore, to complete the proof, we construct $W_0$ in Lemma \ref{W_0_construction} for cases $(1)$ and $(3)$ above, and only for $b\leq \nu$ in case $(2)$ ($m=0$ is not applicable when $b\geq \nu+1$).
\end{proof}

\begin{theorem}{\label{contr_factor}}
Let $r= s+p^t(p-1)d$, with $p\nmid d$, $s =b+c(p-1)$ and suppose also that $2\leq b\leq p, \ 2\leq c\leq p-1$.  Fix $a_p$ such that $\nu(a_p)$ is non-integral, $1<\nu(a_p)<c-\epsilon$ and $t>\nu(a_p)+c$. Assume \textbf{Conjecture \ref{mono II conj}} is true. Then the map $P$ surjects from $\indkg\left(\frac{V_r^{(n_0)}}{V_r^{(n_0 +1)}}\right)$, where $n_0$ is defined as follows:
\begin{enumerate}
\item If $b\geq c+\nu-1$, then $n_0 = \nu$.
\item If $c-1\leq b\leq c+\nu-2$ and $\nu = c-2$, then $n_0 = b-c+1$.
\item  Suppose $2\leq b\leq c-2$ and $1\leq \nu\leq \text{min}\{c-2,p+b-c\}$. Then $n_0 =\nu$ if $b\leq \nu$ and $n_0 = 0$ if $b = \nu+1$.
\item If $2\leq b\leq 2c-4-p$ and $\nu\in\{p+b-c+1, c-3\}$, then $n_0 = p+b-c+1$. 
\end{enumerate}
\end{theorem}


\section{Main result}

In the following proposition, we determine $\vkapL$ for all $k'>k$, where $k'=r+2,\ k=s+2$ and $r,\ s$ are as given in the statement below. We apply Lemma \ref{vrc 1} in the proposition below to consider only those $b, c, \nu$ where the reduction is necessarily irreducible. In particular, we have removed the point $b=2c-4-p$ coming from Theorem \ref{contr_factor} $(3)$ since it is a case with possibly reducible reduction $\vkapL$.


\begin{prop}\label{vk'ap conj}    
Let $p \geq 5$. Let $r= s+p^t(p-1)d$, with $p\nmid d$, $s =b+c(p-1)$ and suppose also that $2\leq b\leq p, \ 2\leq c\leq p-1$.  Fix $a_p$ such that $\nu(a_p)$ is non-integral, $1<\nu(a_p)<c-\epsilon$ and $t>\nu(a_p)+c$. Assume \textbf{Conjecture \ref{mono II conj}} is true and let $k'=r+2$. 
\begin{enumerate}
\item If $b\geq c+\nu-1$ such that $b\not=2\nu+1$ and $(b,\nu) \not= (p,1)$, then $\vkapL\cong\ind\left(\omega^{b+1+\nu(p-1)}_2\right)$.
\item If $c-1\leq b\leq c+\nu-2$ and $\nu = c-2$, then $\vkapL\cong\ind\left(\omega_2^{b+1+(b-c+1)(p-1)}\right)$.
\item  Suppose $2\leq b\leq c-2$ and $1\leq \nu\leq \text{min}\{c-2,p+b-c\}$. Also assume that $b\not=2\nu+2-p$ if $b\leq 2c-2-p$. Then $\vkapL\cong\ind\left(\omega_2^{b+1+(\nu+1)(p-1)}\right)$ if $b\leq \nu$, and $\vkapL\cong\ind\left(\omega_2^{b+1}\right)$ if $b=\nu+1$.
\item If $2\leq b\leq 2c-5-p$ and $\nu\in\{p+b-c+1, c-3\}$, then $\vkapL\cong\ind\left(\omega_2^{k_0}\right)$, where $k_0=b+1+(p+b-c+2)(p-1)$.
\end{enumerate}
\end{prop}

\begin{proof}
We prove the corollary in the following parts.\\
\textbf{Case} $(1)$. $c+\nu-1\leq b\leq p$ \\
By Theorem \ref{contr_factor}, we have 
$$P: \indkg\left(\frac{V^{(\nu)}_r}{V^{(\nu+1)}_r}\right)\twoheadrightarrow\bar{\Theta}_{k',a_p}.$$
Observe that $b\geq 2\nu$ as by assumption $b\geq \nu+c-1$, and $\nu\leq c-1$. We also observe that the equality $b = 2\nu$ occurs only for $\nu=c-1$. But $b= 2c-2$ gives $\epsilon = 1$, and so we must have by hypothesis that $\nu \leq c-2$. Thus, as such we must have $b \geq 2 \nu +1$. Therefore, if $b\not=2\nu+1$ and $(b,\nu) \not= (p,1)$, then by Lemma \ref{vrc 1} (with $n=\nu$), we have 
$$\vkapL\cong \ind\left(\omega_2^{b+1+\nu(p-1)}\right) \ \text{as} \ 2\nu+2 \leq b\leq p $$ 

\textbf{Case} $(2)$. $c-1\leq b\leq c+\nu-2$ and $\nu =c-2$ (so $b\leq 2c-4$)\\
By Theorem \ref{contr_factor}, we have
$$P:\indkg\left(\frac{V_r^{(b-c+1)}}{V_r^{(b-c+2)}}\right)\twoheadrightarrow\bar{\Theta}_{k', a_p}.$$ 
\textbf{Subcase ($i$)} $b=c-1$\\
If $b\not=p-2$, then Lemma \ref{vrc 1} (with $n=0$) gives $\vkapL\cong\ind\left(\omega_2^{b+1}\right)$. If $b=p-2$, then by (\ref{r' not p}) (with $n = 0$ and $r' = 2p-3$), we see that the image of $\indkg\left(V_{p-2}\right)$ in $\indkg\left(\frac{V_r}{V^{(1)}_r}\right)$ is generated by $[1,\ x^r]$ and the latter belongs to $\kerp$ by Remark 4.4 in \cite{KBG}. Hence, $P$ factors through $\ind^G_{KZ}\left(V_1\otimes D^{p-2}\right)$. Therefore, Proposition 3.3 in \cite{KBG} gives $\vkapL\cong \ind\left(\omega^{2+(p-2)(p+1)}_2\right)$.  We conclude by observing that $\omega^{2+(p-2)(p+1)}_2$ is conjugate to $\omega^{b+1}_2$ as $b = p-2$ and $2+(p-2)(p+1)= p(b+1)$. \linebreak\\
\textbf{Subcase ($ii$)} $c\leq b\leq 2c-4$\\
In this case, we note that $b\geq 2n+2$ with $n=b-c+1$ as $b\leq 2c-4$, and also note that $(b,n)\not=(p,1)$ as $c\leq p-1$. Hence, Lemma \ref{vrc 1} (with $n=b-c+1$) gives $\vkapL\cong\ind\left(\omega_2^{b+1+(b-c+1)(p-1)}\right)$.\\ 
\textbf{Case} $(3)$. $2\leq b\leq c-2$ and $1\leq v\leq \text{min}\{c-2,p+b-c\}$\\
\textbf{Subcase ($i$)} $b\leq \nu\leq \text{min}\{c-2,p+b-c\}$\\
By Theorem \ref{contr_factor}, we get 
$$P:\indkg\left(\frac{V_r^{(\nu)}}{V_r^{(\nu+1)}}\right)\twoheadrightarrow\bar{\Theta}_{k',a_p}.$$
Observe that $2\nu+1-(p-1)\leq b\leq 2\nu$ as $b\leq \nu\leq\text{min}\{c-2, p+b-c\}$. The equality $b=2\nu+2-p$ occurs only if $b\leq 2c-2-p$, and these are possibly reducible cases. Therefore, by Lemma \ref{vrc 1}, we have $\vkapL\cong\ind\left(\omega_2^{b+1+(\nu+1)(p-1)}\right)$ as $2\nu+1-(p-1)<b\leq 2\nu$.\\
\textbf{Subcase ($ii$)} $b =\nu+1$ and $1\leq \nu\leq \text{min}\{c-2,p+b-c\}$\\
By Theorem \ref{contr_factor} we get
$$P:\indkg\left(\frac{V_r}{V_r^{(1)}}\right)\twoheadrightarrow\bar{\Theta}_{k', a_p}.$$ 
If $b\not=p-2$, then Lemma \ref{vrc 1} gives $\vkapL\cong\ind\left(\omega_2^{b+1}\right)$. Observe that $b=p-2$ occurs only if $(\nu,c)=(c-2,p-1)$, and in this case, we proceed exactly the same as in the proof of part ($1$) above to get the required result.\\ 
\textbf{Case} $(4)$. $2\leq b\leq2c-5-p$ and $\nu\in\{p+b-c+1, c-3\}$\\
In this case, Theorem \ref{contr_factor} gives
$$P:\indkg\left(\frac{V_r^{(p+b-c+1)}}{V_r^{(p+b-c+2)}}\right)\twoheadrightarrow\bar{\Theta}_{k',a_p}.$$
Observe that if $n=p+b-c+1$, then $2n+1-(p-1)<b<2n$ as $b\leq 2c-5-p$. Hence, using Lemma \ref{vrc 1}, we obtain $\vkapL\cong\ind\left(\omega_2^{b+1+(p+b-c+2)(p-1)}\right)$. 
\end{proof}
The following theorem is the main result of this section, where we prove local constancy by assuming Conjecture \ref{mono II conj}. 
\begin{theorem}\label{main_result III}
Let $k = b+c(p-1)+2$ with $2\leq c\leq p-1$, $2\leq b\leq p$ and $p\geq 5$. Fix $a_p$ such that $\nu(a_p)$ is non-integral, $1<\nu(a_p)<c-\epsilon$, and let $t>\nu(a_p)+c$. Assume \textbf{Conjecture \ref{mono II conj}} is true. 
\begin{enumerate}
\item If $b\geq c+\nu-1$ such that $b\not=2\nu+1$ and $(b,\nu)\not=(p,1)$, then $\vkapL\cong\ind\left(\omega^{b+1+ \nu(p-1)}_2\right)$ for all $k'\in k+p^t(p-1)\mathbb{Z}_{\geq 0}$.
\item If $c-1\leq b\leq 2c-4$ and $\nu = c-2$, then $\vkapL\cong\ind\left(\omega_2^{b+1+(b-c+1)(p-1)}\right)$ for all $k'\in k+p^t(p-1)\mathbb{Z}_{\geq 0}$.
\item  Suppose $2\leq b\leq c-2$ and $1\leq \nu\leq \text{min}\{c-2,p+b-c\}$. Also assume that $b\not=2\nu+2-p$ if $b\leq 2c-2-p$. Then for all $k'\in k+p^t(p-1)\mathbb{Z}_{\geq 0}$, $\vkapL\cong\ind\left(\omega_2^{b+1+(\nu+1)(p-1)}\right)$ if $b\leq \nu$, and $\vkapL\cong\ind\left(\omega_2^{b+1}\right)$ if $b=\nu+1$.
\item  If $2\leq b\leq 2c-5-p$ and $\nu\in\{p+b-c+1, c-3\}$, then $\vkapL\cong\ind\left(\omega_2^{k_0}\right)$ for all $k'\in k+p^t(p-1)\mathbb{Z}_{\geq 0}$, where $k_0=b+1+(p+b-c+2)(p-1)$.
\end{enumerate} 
In the above cases the Berger's constant $m(k,a_p)$ exists and $m(k,a_p)\leq \lceil\nu(a_p)\rceil+c+1$. 
\end{theorem} 
\begin{rk}\label{rk_uncodi_1}
    Here, we make an important remark that statement $(1)$ is unconditional on Conjecture \ref{mono II conj}. This is because Conjecture \ref{mono II conj} is a consequence of Conjecture \ref{mat_conj}, and the latter is proved in Lemma \ref{proof_conj_case1} for the cases required in statement $(1)$ above (see also Remark \ref{conj_required_1} above).\
\end{rk}
\begin{proof}
Observe that since $\nu(a_p)< c-\epsilon$ we have
\begin{align*}
3\nu(a_p)+\frac{(k-1)p}{(p-1)^2}+1 &<4c+1-3\epsilon+\frac{(b+1)p}{(p-1)^2}+\frac{c}{(p-1)}\\
&<b+c(p-1)+2 =k \quad \forall\ p\geq 5.
\end{align*}
The last inequality follows as $(b+1)(p^2-3p+1)+c(p^2-6p+4)(p-1)+3\epsilon(p-1)^2>0$. Therefore, we get $k>3\nu(a_p)+\frac{(k-1)p}{(p-1)^2}+1$ so, by Theorem \ref{Berger} there exists a constant $m = m(k, a_p)$ such that for all $k'' \in k+p^{m-1}(p-1)\mathbb{Z}_{\geq 0}$, we have $\bar{V}_{k'',a_p}\cong \vkap$. Observe that Proposition \ref{vk'ap conj} determines $\vkapL$ for all $k'\in k+p^{t}(p-1)\mathbb{Z}_{>0}$ (punctured disk around $k = s+2$) with $t>\nu(a_p)+c$ and shows that $\vkapL$ is constant function of $k'$ for each $k$. Note that both the disks $k+p^{m-1}(p-1)\mathbb{Z}_{\geq 0}$ and $k+p^{t}(p-1)\mathbb{Z}_{>0}$ have the same center. Therefore, we have $\vkapL\cong \vkap$ for all $k'\in k+p^{t}(p-1)\mathbb{Z}_{\geq 0}$. Hence, we get $m(k,a_p)\geq t+1$ and also determine $\vkap$ in each case of the theorem. Thus, we have the desired result.
\end{proof}

\begin{cor}\label{new_k}
Fix $a_p$ such that $\nu(a_p)$ is non-integral and $1<\nu(a_p)<p-1$. Let $p\geq 13$  and assume \textbf{Conjecture \ref{mono II conj}} is true. Suppose that $k\not\in\{2\nu+3+c(p-1),2\nu+4-p+c(p-1)\}$ and $(k,\nu)\not= (p+2+c(p-1),1)$, where $c$ is given below. Theorem \ref{main_result III} computes $\vkap$ in the following cases:
\begin{enumerate}
     \item If $\nu = 1$, then for all $$k\in\bigcup\limits_{c = 2+\epsilon}^{p-1}[c+2+c(p-1), p+2+c(p-1)]\bigcup\limits_{c = 3}^{p-1}\{4+c(p-1)\}.$$
    \item If $\nu = 2$, then for all $$k\in \bigcup\limits_{c= 3+\epsilon}^{p-1}[c+3+c(p-1),p+2+c(p-1)]\bigcup\limits_{c=4}^{p-1}[4+c(p-1),5+c(p-1)]\bigcup\{4p+2\}.$$
    \item If $\nu = 3$, then for all $$k\in \bigcup\limits_{c= 4+\epsilon}^{p-2}[c+4+c(p-1),p+2+c(p-1)]\bigcup\limits_{c=5}^{p-1}[4+c(p-1),6+c(p-1)]\bigcup\{5p+2,5p+3\}.$$
    \item If $4\leq \nu \leq\frac{p-1}{2}$, then for all 
    \begin{align*}
        &k\in \bigcup\limits_{c= \nu+1+\epsilon}^{p+1-\nu}[c+\nu+1+c(p-1),p+2+c(p-1)]\bigcup\limits_{c=\nu+2}^{p-\nu+2}[4+c(p-1),\nu+3+c(p-1)]\\ &\bigcup\limits_{c = p-\nu+3}^{p-1}[c+\nu+1-p+c(p-1), \nu+3+c(p-1)]\bigcup[2+(\nu+2)p, \nu+(\nu+2)p].
    \end{align*}
    \item If $\frac{p+1}{2}\leq \nu\leq p-7$, then for all $$\hspace{2em}k\in\bigcup\limits_{c= \nu+3}^{p-1}[c+\nu+1-p+c(p-1), \nu+3+c(p-1)]\bigcup[(\nu+1)p+\nu+2, (\nu+1)p+\nu+p].$$
    \item If $p-6\leq \nu\leq p-5$, then for all  
    \begin{align*}
       \hspace{2em}&k\in\bigcup\limits_{c = \nu+2}^{\frac{p+\nu+3}{2}}[2c-p+c(p-1),\nu+3+c(p-1)]\bigcup\limits_{c= \nu+4}^{p-1}[c+\nu+1-p+c(p-1),2c-1-p+c(p-1)]\\&\hspace{2em}\bigcup[\nu+4+(\nu+2)(p-1),2\nu+2+(\nu+2)(p-1)]\bigcup\{2\nu+4+(\nu+2)(p-1)\}.
    \end{align*}
    \item If $\nu = p-4$, then for all $$\hspace{2em}k\in\bigcup\limits_{c = p-2}^{p-1}[2c-p+c(p-1), p-1+c(p-1)]\bigcup[1+(p-1)^2, p-5+(p-1)^2]\bigcup\{p-3+(p-1)^2\}.$$
    \item If $\nu = p-3$, then for all $k\in[p-2+(p-1)^2, p+2+(p-1)^2]$
\end{enumerate}     
\end{cor}
\begin{rk} We note that $\nu = p-2$ cannot occur. This is because we must also have $\nu = p-2 \leq c-1-\epsilon$, and so $\epsilon = 0$. Thus, $\nu = p-2$ could come only from Theorem \ref{main_result III} $(1)$ wherein $p \geq b \geq c+\nu -1$, forcing $p \leq 4$.
    
\end{rk}
\section{Proof of more cases of Conjecture \ref{mat_conj}}\label{pf_more_cases_conj}

In proving the following lemma, the main observation is to obtain some simple expression for $\alpha(i,l)$, which is essentially equivalent to giving some simple expression for $\alpha_1(i,l)$ (see \eqref{dfn alpha}). Lemma \ref{alpha cong}
gives that $\alpha_1(i,l)\equiv\alpha'(i,l)$ (for definition of $\alpha'(i,l)$ see \eqref{dfn alph'}). To prove $(1)$ - $(3)$ of the theorem below, we show modulo $p$, $\alpha'(i,l)= 0$ if $l\geq i$ and $\alpha'(i,l)=(-1)^{b-m+l}u_lv_i(m+1-i){{i-1}\choose l}$ if $l\leq i-1$ by using Lemmas \ref{alpha' cong app} \& \ref{mat cong}. To prove the last part of the following lemma, we use Lemmas \ref{alpha' cong II} \& \ref{alpha' cong II app} to do a similar analysis. 
Recall that $m\in[1, c-1-\epsilon]$ in Conjecture \ref{mat_conj}.
\begin{theorem}\label{proof_conj}
    Conjecture \ref{mat_conj} is true in the following cases (with the additional condition on $m$):
    \begin{enumerate}
        \item $m\in[1, b-c]$ if $b\geq c+1$.
        \item $m\in[1, \frac{b-1}{2}]\cap[b-c+1, c-1-\epsilon]$ if $b\geq c$.
        \item $m\in\left([1, \frac{b-1}{2}]\cup[b,c-2]\right)\cap[1, p+b-c]$ if $2\leq b\leq c-1$.
        \item $m\in[p+b-c+1, \frac{p+b-2}{2}]\cap[p+b-c+1, c-3]$ if $2\leq b\leq 2c-3-p$.
    \end{enumerate}
\end{theorem}
Before we give the proof of the theorem above, we would like to state a corollary of it and make some remarks.
\begin{cor}\label{thm_conj}
Conjecture \ref{mat_conj} is true in the following cases (with the additional condition on $m$):
\begin{enumerate}
\item For all $m$ if $b\geq 2c-3$.
\item $m\in[1,\frac{b-1}{2}]\cup[b,c-2]$ if $2c-2-p\leq b\leq 2c-4$.
\item $m\in[1, \frac{b-1}{2}]\cup[b, c-3]$ if $2c-4-p\leq b\leq 2c-3-p$.
\item $m\in[1, \frac{b-1}{2}]\cup[b, \frac{p+b-2}{2}]$ if $2\leq b\leq 2c-5-p$.
\end{enumerate}
\end{cor}
\begin{proof}
When $b\geq 2c-2$, the full range of $m$ i.e, $[1,c-1-\epsilon]$ is contained in $[1,b-c]$ as $c-1-\epsilon\leq b-c$, and so Theorem \ref{proof_conj} $(1)$ gives Corollary \ref{thm_conj} $(1)$ in this case. For $b = 2c-3$ we use Theorem \ref{proof_conj} $(1)$ $ \& $ $(2)$ for $1\leq m\leq c-3$ and $m = c-2$ respectively to prove Corollary \ref{thm_conj} $(1)$. Next, when $c\leq b\leq 2c-4$, observe that $[b,c-2]$ is empty and $[1,\frac{b-1}{2}] = [1,b-c]\cup([1,\frac{b-1}{2}]\cap[b-c+1,c-1-\epsilon])$ as $b-c< \frac{b-1}{2}$ (since $b\leq 2c-4$). Therefore, the statement $(2)$ of Corollary \ref{thm_conj} is obtained by Theorem \ref{proof_conj} $(1)$ $\&$\ $(2)$ when $c\leq b\leq 2c-4$ and only by Theorem \ref{proof_conj} $(3)$ when $2c-2-p\leq b\leq c-1$. Lastly, note that $\left([1, \frac{b-1}{2}]\cup[b,c-2]\right)\cap[1, p+b-c] = [1,\frac{b-1}{2}]\cup[b,p+b-c]$ if $b\leq 2c-3-p$. Thus, the last two statements of the corollary follow by Theorem \ref{proof_conj} $(3)$ \& $(4)$ as $[p+b-c+1, \frac{p+b-2}{2}]\cap[p+b-c+1, c-3]$ is equal to $[p+b-c+1, c-3]$ if $b\geq 2c-4-p$ and equal to $[p+b-c+1, \frac{p+b-2}{2}]$ if $b\leq 2c-5-p$.
\end{proof}
\begin{rk}\label{remaining_conjecture}
From the above corollary, we see that Conjecture \ref{mat_conj} remains to be proven in the following cases:
\begin{enumerate}
\item $m\in[\frac{b}{2}, b-1]\cap[1,c-2]$ if $2c-2-p\leq b\leq 2c-4$.
\item $m\in[\frac{b}{2}, b-1]$ if $2c-4-p\leq b\leq 2c-3-p$.
\item $m\in[\frac{b}{2}, b-1]\cup[\frac{p+b-1}{2}, c-3]$ if $2\leq b\leq 2c-5-p$.
\end{enumerate}
Using SageMath, we have verified that the matrix $A$ defined in \eqref{matrix_A_conj} is invertible for the above ranges of $b,c,m$ and primes $p$ up to $97$.
\end{rk}
We now prove Theorem \ref{proof_conj}. It will be clear from the proof below that the matrix $A$ defined in \eqref{matrix_A_conj} is non-invertible for the ranges of $b,c,m$ in Theorem \ref{proof_conj}, and moreover, the last row of $A$ is zero $\bmod\ p$.
\begin{proof}
Recall that from \S \ref{conje_local_sec}, we have the matrix 
   $$A=\left(\alpha(i,l)\right)_{\substack{1\leq i\leq m+1\\ 0\leq l\leq m}}$$
   where $\alpha(i,l)$ is as defined in \eqref{dfn alpha}. Now, we show that the system of linear equations $AX=e_i(\bmod p)$ has a solution for all $1\leq i\leq m$.\\
   \textbf{Case (i)} $b\geq c+1$ and $m\in[1, b-c]$. This is proved in Lemma \ref{proof_conj_case1}.\\
\textbf{Case (ii)}  $b\geq c$ and $m\in[1, \frac{b-1}{2}]\cap[b-c+1, c-1-\epsilon]$\\
    In this case, we write $A$ as follows
	$$A= \begin{pmatrix}
		  A' & B'\\
		  A'' & B''
		  \end{pmatrix}$$
	where the ranges of $i$ and $l$ are divided into non-empty intervals $[1,\ b-c+1], \ [b-c+2, \ m+1]$ and $[0, b-c],\ [b-c+1,\ m]$ respectively, determining the order of blocks.\\
	\textbf{Subcase (a)} $0\leq l\leq b-c$ and $1\leq i\leq b-c+1$\\
	For the above ranges of $i$ and $l$, Lemma \ref{alpha cong} gives $\alpha_1(i,l)\equiv0\bmod p$ and 
    $${{r-l}\choose b-m+(i+c-m-1)(p-1)}\equiv {{b-c-l}\choose b-c+1-i}{c\choose i+c-m-1}.$$ 
    Hence, we have
	\begin{eqnarray*}
	   \alpha(i,l)&\equiv & {{b-c-l}\choose b-c+1-i}{c\choose c-m-1+i}\\
				  &\equiv & 0\bmod p\iff i<l+1.
	\end{eqnarray*}
    According to the calculation above, modulo $p$ the matrix $A'$ is an invertible lower triangular matrix. \\
    \textbf{Subcase (b)} $b-c+1\leq l\leq m$ and $b-c+2\leq i\leq m+1$\\
	In this case, we note that $m\leq p+b-c$ as $b\geq c$ and $m\leq c-1-\epsilon\leq p-1$. Using Lemma \ref{alpha cong}, we get 
	$$\alpha(i,l)\equiv \begin{cases}
	\alpha'(i, l)+{{p+b-c-l}\choose p+b-c+1-i}{{c-1}\choose c-m-2+i} & \text{if}\ b-c+1\leq l\leq m,\ b-c+2\leq i\leq m\\
	 \alpha'(i, l) & \text{if}\ b-c+1\leq l\leq m, i = m+1
	\end{cases}$$
	where $\alpha'(i,l)$ is as defined in \eqref{dfn alph'}. By Lemma \ref{alpha' cong app}, we have 
		$$\alpha'(i,l) \equiv u_lv_i\underset{\epsilon_1\leq k\leq c_1}\sum (-1)^k(c-k){{c-m-1+i}\choose k}{{i+c-m-2-k}\choose i-(b-c+2)}{{l+c-m-1-k}\choose l},$$                    where $c_1=\text{min}\{c-m-1, b-m\}$, $u_l = (-1)^l l!(p+b-c-l)!(c-1)!$ and
		  \begin{eqnarray*}
			v_i &=& \frac{(-1)^{i+1+c-m}(i-(b-c+2))!}{(i-1)! (m+1-i)! (c-m-1+i)!}
			\end{eqnarray*}
            We note that $b\geq 2m+1$. Using the above calculations and Lemma \ref{mat cong} along with Remark \ref{rk of mat cong}, we obtain 
            $$\alpha'(i,l)\equiv\begin{cases}
                0 &\text{if}\ i\leq l\\
                (-1)^{b-m+l}u_lv_i(m+1-i){{i-1}\choose l} & \text{if}\ l\leq i-1.
            \end{cases}$$
			Hence, we get for $b-c+1 \leq l \leq m$
			\begin{align*}
			    \hspace*{2em}\alpha(i,l)\equiv \begin{cases}
			{{p+b-c-l}\choose p+b-c+1-i}{{c-1}\choose c-m-2+i}\quad  \text{if}\ b-c+2\leq i\leq m,\ i\leq l\\
			(-1)^{b-m+l}u_lv_i(m+1-i){{i-1}\choose l}+{{p+b-c-l}\choose p+b-c+1-i}{{c-1}\choose c-m-2+i}\ \ \text{if}\ c_1'\leq i\leq m,\  l\leq i-1\\
			0 \quad \text{if}\ i = m+1
			\end{cases}
			\end{align*}
            where $c_1'=b-c+2$. The last case is clear as $(-1)^{b-m+l}u_lv_i(m+1-i){{i-1}\choose l}=0$ for $i=m+1$. For $i\leq m$,  we note that 
						\begin{align*}
							&	{{p+b-c-l}\choose p+b-c+1-i}{{c-1}\choose c-m-2+i}\\
							&= \frac{(-1)^{b+c+1+i}l!(p+b-c-l)!(c-1)!(i-(b-c+2))!{{i-1}\choose l}}{(i-1)!(m+1-i)!(c-m-2+i)!}.
						\end{align*}         
						Thus, we have for $b-c+1 \leq l \leq m$              
						\begin{align*}
							& \alpha(i,l)\equiv \begin{cases}
								\frac{(-1)^{b+c+1+i}l!(p+b-c-l)!(c-1)!(i-(b-c+2))!}{(i-1)!(m+1-i)!(c-m-2+i)!}{{i-1}\choose l} & \text{if}\ b-c+2\leq i\leq m,\ i\leq l\vspace*{2mm}\\
								\frac{(-1)^{b+c+i+1}c\ l!(p+b-c-l)!(c-1)!(i-(b-c+2))!}{(i-1)!(m+1-i)!(c-m-1+i)!}{{i-1}\choose l} &\text{if}\ b-c+2\leq i\leq m,\  l\leq i-1\\
								0 & \text{if}\ i = m+1
							\end{cases}
						\end{align*}
						\begin{align*}
							&\equiv  \begin{cases}
								u'_lv'_i{{i-1}\choose l} & \text{if}\ b-c+2\leq i\leq m,\ b-c+1\leq l\leq m \vspace*{2mm}\\
								0 & \text{if}\ i = m+1
							\end{cases}
						\end{align*}
						where $u'_l = l!(p+b-c-l)!$ and 
						$$v'_i = \begin{cases}
							\frac{(-1)^{b+c+1+i}(c-1)!(i-(b-c+2))!}{(i-1)!(m+1-i)!(c-m-2+i)!} & \text{if}\ b-c+2\leq i\leq m,\ i\leq l\vspace*{2mm}\\
							\frac{(-1)^{b+c+i+1}c\ (c-1)!(i-(b-c+2))!}{(i-1)!(m+1-i)!(c-m-1+i)!} &\text{if}\ b-c+2\leq i\leq m,\  l\leq i-1.\\
						\end{cases}$$
						From the above congruence on $\alpha(i,l)$, we have 
						$$B''\equiv\begin{pmatrix}
							B''' & 0\\
							0 & 0
						\end{pmatrix}\bmod p, \quad \text{where}\ B''' = \left(u'_lv'_i{{i-1}\choose l}\right)_{\substack{b-c+2\leq i\leq m\\ b-c+1\leq l\leq m-1}}$$
					   and $0$ denotes the zero matrix of the required size.
						We note that $B'''$ is invertible (lower triangular with non-zero diagonal entries).
      
We finally also note that by Lemma \ref{alpha cong} the matrix $A''\equiv 0\bmod p$ since both $\alpha'(i,l)$ and $X(i,l)$ are zero $\bmod \ p$ when $b-c+2 \leq i \leq m+1$ and $0\leq l \leq b-c$. Thus, 
						$$A\equiv\begin{pmatrix}
							A' & B'\\
							0 & B''
						\end{pmatrix} $$    
						where $A'$ is invertible $\bmod p$, and $B''$ is as above. Hence for every $1\leq i'\leq m$, modulo $p$ the row rank  of $[A: \textbf{e}_{i'}]$ is the same as the row rank of $A$ (which is $m$). Thus, the linear systems $AX=e_{i'}(\bmod p)$ has a solution for all $1\leq i'\leq m$. \\
                        \textbf{Case (iii)} $2\leq b\leq c-1$ and $m\in\left([1, \frac{b-1}{2}]\cup[b,c-2]\right)\cap[1, p+b-c]$\vspace{.3em}\\
                        Observe that $i\in[1,m+1]\subset[b-c+2, p+b-c+1]$ and $l\in[0, m]\subset[b-c+1, p+b-c]$. Therefore, Lemma \ref{alpha cong} gives
                        $$\alpha(i,l)\equiv \begin{cases}
							\alpha'(i, l)+{{p+b-c-l}\choose p+b-c+1-i}{{c-1}\choose c-m-2+i} & \text{if}\ 0\leq l\leq m,\ 1\leq i\leq m\\
							\alpha'(i, l) & \text{if}\ 0\leq l\leq m, i = m+1
						\end{cases}$$
where $\alpha'(i,l)$ is as defined in \eqref{dfn alph'}. If $m\in[1,\frac{b-1}{2}]$, then all the computations of Subcase (b) of Case (ii) for $\alpha(i,l)$ will carry over here and give that the matrix $A$ is of the form $B''$ i.e.,
$$A\equiv\begin{pmatrix}
B''' & 0\\
0 & 0
\end{pmatrix}\bmod p, \quad \text{where}\ B''' = \left(u'_lv'_i{{i-1}\choose l}\right)_{\substack{1\leq i\leq m\\ 0\leq l\leq m-1}}.$$
Thus, following the same argument as in the previous case, we get that the linear systems $AX=e_{i'}(\bmod p)$ has a solution for all $1\leq i'\leq m$. \\
If $m\in[b, c-2]$, then using Lemma \ref{alpha cong}, we observe that $\alpha_1(i,l)\equiv 0\bmod p$ (indeed in this case, the summation in \eqref{dfn alph'} of $\alpha'(i,l)$ is empty as $\epsilon_1 \geq 1$). Hence, we have
$$\alpha(i,l)\equiv 
\begin{cases}
{{p+b-c-l}\choose p+b-c+1-i}{{c-1}\choose c-m-2+i} & \text{if}\ 0\leq l\leq m,\ 1\leq i\leq m\\
0 & \text{if}\ 0\leq l\leq m, i = m+1.
\end{cases}$$
For the above values of $l$ and $i$, we note that ${{p+b-c-l}\choose p+b-c+1-i}=0\ \ \iff\ i<l+1$. Therefore, we can see $A$ as a block matrix of the form given in \eqref{pf_conj_case_1_eqn}, where modulo $p$, $A'$ is invertible and $A'', B''$ are zero. Thus, the linear systems $AX=e_{i'}(\bmod p)$ has a solution for all $1\leq i'\leq m$.\\
\textbf{Case (iv)} $2\leq b\leq 2c-3-p$ and $m\in[p+b-c+1, \frac{p+b-2}{2}]\cap[p+b-c+1, c-3]$\\
In this case, we note that $\epsilon_1= 1$ or $2$ as $b\leq m$. If $m\geq p+b-c+2$, then we write $A$ as follows
	$$A= \begin{pmatrix}
		  A' & B'\\
		  A'' & B''\\
          A''' & B'''
		  \end{pmatrix}$$
	where the ranges of $i$ and $l$ are divided into non-empty intervals $[1,\ p+b-c+1], \ [p+b-c+2, \ m],\ [m+1,m+1]$ and $[0, p+b-c],\ [p+b-c+1,\ m]$ respectively, determining the order of blocks.\\
    If $m=p+b-c+1$, then we write $A$ as follows
    $$A= \begin{pmatrix}
		  A' & B'\\
          A''' & B'''
		  \end{pmatrix}$$
	where the ranges of $i$ and $l$ are divided into non-empty intervals $[1,\ m],\ [m+1,m+1]$ and $[0, m-1],\ [m,\ m]$ respectively, determining the order of blocks.\\
    \textbf{Subcase (a)} $0\leq l\leq p+b-c$\\
     Observe that using Lemma \ref{alpha cong}, we get
     $${{r-l}\choose b-m+(c-m-1+i)(p-1)}\equiv \begin{cases}
          {{p+b-c-l}\choose p+b-c+1-i}{{c-1}\choose c-m-2+i}\quad\text{if}\ \ 1\leq i\leq p+b-c+1\\
          0 \quad\text{if}\ \ p+b-c+2\leq i\leq m+1
      \end{cases}$$
      and $\alpha_1(i,l)\equiv\alpha'(i,l)=0$ (indeed in this case, the summation in \eqref{dfn alph'} of $\alpha'(i,l)$ is empty as $c-m-\epsilon_1<c-b$). Thus, we have
      $$\alpha(i,l)\equiv \begin{cases}
          {{p+b-c-l}\choose p+b-c+1-i}{{c-1}\choose c-m-2+i} &\text{if}\ \ 1\leq i\leq p+b-c+1\\
          0 &\text{if}\ \ p+b-c+2\leq i\leq m+1
      \end{cases}$$

      Note that for $1\leq i\leq p+b-c+1$, we have
      $${{p+b-c-l}\choose p+b-c+1-i}=0\iff\ i<l+1.$$
      Hence, modulo $p$ we get $A'$ is lower triangular invertible (because all the diagonal entries are non-zero modulo $p$, given by $i = l+1$), and $A'',\ A'''$ are zero modulo $p$.\\
      \textbf{Subcase (b)} $p+b-c+1\leq l\leq m$ and $p+b-c+2\leq i\leq m+1$\\
      Here, $ l \in [p+b-c+1, m]$ and $i \in [p+b-c+2, m+1]$. So using \eqref{dfn alph'} and Lemma \ref{alpha' cong II app}, we get 
      \begin{center}
          $\alpha'(i,l)\equiv u_lv_i\underset{\epsilon_1-1\leq k\leq c-m-2}\sum(-1)^k(c-1-k){{c-m-2+i}\choose k}{{i+c-m-3-k}\choose i-(p+b-c+2)}{{l+c-m-2-k}\choose l},$
      \end{center}
      where $u_l=(-1)^ll!(2p+b-c-l)!(c-2)!$ and
      $$v_i=\frac{(-1)^{i+c-m}(i-(p+b-c+2))!}{(i-1)!(m+1-i)!(c-m-2+i)!}.$$
      Since $m\leq \frac{p+b-2}{2}$ it follows that $b\geq 2m+1-(p-1)$ and $\epsilon_1 = 1$. Using the above calculations and Lemma \ref{alpha' cong II} along with Remark \ref{rk of alpha' cong II}, we obtain  
              $$\alpha'(i,l)\equiv\begin{cases}
                  0 & \text{if}\ l\geq i\\
                  (-1)^{b-m+l}u_lv_i(m+1-i){{i-1}\choose l}& \text{if}\ l\leq i-1.
              \end{cases}$$
      Hence, by using Lemma \ref{alpha cong}, we have
      \begin{align*}
          &\alpha(i,l)\equiv \begin{cases}
          {{2p+b-c-l}\choose 2p+b-c+1-i}{{c-2}\choose c-m-3+i}\quad\text{if}\ \ p+b-c+2\leq i\leq m,\ i\leq l\\
          (-1)^{b-m+l}u_lv_i(m+1-i){{i-1}\choose l}+{{2p+b-c-l}\choose 2p+b-c+1-i}{{c-2}\choose c-m-3+i}\ \text{if}\ \alpha_0\leq i\leq m,\ l\leq i-1 \\
          0 \quad \text{if}\ i=m+1,
      \end{cases}
      \end{align*}
      where $\alpha_0=p+b-c+2$. \\
      Note that for $p+b-c+2\leq i\leq m$, we have
      \begin{align*}
          &{{2p+b-c-l}\choose 2p+b-c+1-i}{{c-2}\choose c-m-3+i}\\
          &\equiv\frac{(-1)^{b-c+i}l!(2p+b-c-l)!(i-(p+b-c+2))!(c-2)!{{i-1}\choose l}}{(i-1)!(m+1-i)!(c-m-3+i)!}
          \end{align*}
          Therefore, we get
          $$\alpha(i,l)\equiv \begin{cases}
          u'_lv'_i{{i-1}\choose l}\quad\text{if}\ \ p+b-c+2\leq i\leq m\\
          0 \quad \text{if}\ i=m+1,
      \end{cases}$$
       where $u_l'=l!(2p+b-c-l)!$ and 
       $$v'_l=\begin{cases}
           \frac{(-1)^{b-c+i}(i-(p+b-c+2))!(c-2)!}{(m+1-i)!(c-m-3+i)!(i-1)!}\ \ \text{if}\ i\leq l\vspace{.4em}\\
           \frac{(-1)^{b-c+i}(i-(p+b-c+2))!(c-1)!}{(m+1-i)!(c-m-2+i)!(i-1)!}\ \ \text{if}\ l\leq i-1.
       \end{cases}$$
       Hence, $B'''$ is zero modulo $p$ and
       $$B''\equiv\begin{pmatrix}
							B''_1 & 0
						\end{pmatrix}\bmod p, \quad \text{where}\ B''_1 = \left(u'_lv'_i{{i-1}\choose l}\right)_{\substack{p+b-c+2\leq i\leq m\\ p+b-c+1\leq l\leq m-1}}.$$
      Hence, modulo $p$ the matrix $B''_1$ is lower triangular invertible.\\

     Using the above calculations, we get 
     $$A\equiv \begin{pmatrix}
		  A' & B'\\
		  0& B''\\
          0 & 0
		  \end{pmatrix},$$
    where modulo $p$ the matrix $A'$ is an invertible, and $B''$ is full rank ($A'''\equiv0\bmod\ p$ and $B'''\equiv 0\bmod\ p$ are just row matrix). Hence, for every $1\leq i'\leq m$, modulo $p$ the row rank  of $[A: \textbf{e}_{i'}]$ is same as the row rank of $A$. Thus, the linear systems $AX=e_{i'}(\bmod p)$ has a solution for all $1\leq i'\leq m$. 
\end{proof}

\noindent\textbf{Acknowledgements.}
We owe a great debt to the work in \cite{maam} and \cite{Berger12} critical to our work. The authors would like to express sincere gratitude to Shalini Bhattacharya for giving useful suggestions regarding this problem. The second author acknowledges the valuable support from IISER Tirupati during his postdoctoral research fellowship at the institute.


 \appendix
 \section{}

\begin{lemma}\label{alpha' cong app}
 Let $2\leq b\leq p$, $0\leq c\leq p-1$ and $1\leq m\leq c-1$. For $\text{max}\{1, b-c+2\}\leq i\leq m+1$ and $\text{max}\{0, b-c+1\}\leq l\leq \text{min}\{m, p+b-c\}$ then
    $$\alpha'(i,l) \equiv u_lv_i\underset{\epsilon_1\leq k\leq c_1}\sum (-1)^k(c-k){{c-m-1+i}\choose k}{{i+c-m-2-k}\choose i-(b-c+2)}{{l+c-m-1-k}\choose l},$$                
where $\alpha'(i,l)$ defined in \eqref{dfn alph'}, $c_1=\text{min}\{c-m-1, b-m\}$, $u_l= (-1)^l l!(p+b-c-l)!(c-1)!$ and \begin{center}
    $v_i =\frac{(-1)^{i+1+c-m}(i-(b-c+2))!}{(i-1)! (m+1-i)! (c-m-1+i)!}.$
\end{center}
\end{lemma}
\begin{proof}
Recall from the definition of $\alpha'(i,l)$ in \eqref{dfn alph'}, we have
$$\alpha'(i,l) = (-1)^{i+1}\underset{\text{max}\{c-b,1\}\leq a\leq c-m-\epsilon_1}\sum{{p+b-c-l}\choose b-c+a}{{c-1}\choose c-m-a}\beta(a,i)$$
as $b-c+1\leq l\leq p+b-c$. Using Lemma \ref{beta_equivalent}, we get
$$\alpha'(i,l) = (-1)^{i+1}\underset{\text{max}\{c-b,1\}\leq a\leq c-m-\epsilon_1}\sum{{p+b-c-l}\choose b-c+a}{{c-1}\choose c-m-a}{{i+a-2}\choose a-1}{{m+a}\choose i+a-1}.$$
Next,
\begin{align*}
    &{{p+b-c-l}\choose b-c+a}{{c-1}\choose c-m-a}{{i+a-2}\choose a-1}{{m+a}\choose i+a-1}\\
    =&\frac{(p+b-c-l)!(c-1)!(i+a-2)!(m+a)!}{(b-c+a)!(p-(l+a))!(c-m-a)!(m+a-1)!(a-1)!(i-1)!(i+a-1)!(m+1-i)!}\\
    \equiv&\frac{(p+b-c-l)!(c-1)!(i+a-2)!(l+a-1)!(m+a)(-1)^{l+a}}{(b-c+a)!(c-m-a)!(a-1)!(i-1)!(i+a-1)!(m+1-i)!}\\
\end{align*}
The last congruence follows by noting that $(p-(l+a))!\equiv\frac{(-1)^{l+a}}{(l+a-1)!}(\bmod\ p)$. On multiplying and dividing by $l!(i+c-m-1)!(i-(b-c+2))!$ on the right-hand side of the above equation, we get
\begin{align*}
    &{{p+b-c-l}\choose b-c+a}{{c-1}\choose c-m-a}{{i+a-2}\choose a-1}{{m+a}\choose i+a-1}\\
    \equiv&(-1)^{i+1+c-m+a}u_lv_i(m+a){{i+c-m-1}\choose c-m-a}{{i+a-2}\choose i-(b-c+2)}{{l+a-1}\choose l}.
\end{align*}
Therefore, we have
\begin{align*}
    \alpha'(i,l)&\equiv u_lv_i\underset{\text{max}\{c-b,1\}\leq a\leq c-m-\epsilon_1}\sum(-1)^{c-m+a} (m+a){{i+c-m-1}\choose c-m-a}{{i+a-2}\choose i-(b-c+2)}{{l+a-1}\choose l}\\
    &=u_lv_i\underset{\epsilon_1\leq k\leq c_1}\sum(-1)^k(c-k){{i+c-m-1}\choose k}{{i+c-m-2-k}\choose i-(b-c+2)}{{l+c-m-1-k}\choose l}.
\end{align*}
The last equality follows by putting $k = c-m-a$.
\end{proof}

Following similar steps as in the proof of Lemma \ref{alpha' cong app}, we obtain the following lemma.
\begin{lemma}\label{alpha' cong II app}
    Let $2\leq b\leq p$, $0\leq c\leq p-1$ and $1\leq m\leq c-2$. For $p+b-c+2\leq i\leq m+1$ and $p+b-c+1\leq l\leq \text{min}\{m, 2p+b-c\}$ then
    \begin{eqnarray*}
        &\alpha'(i,l) \equiv u_lv_i\underset{\epsilon_1-1\leq k\leq c_1}\sum (-1)^k(c-1-k){{c-m-2+i}\choose k}{{i+c-m-3-k}\choose i-(p+b-c+2)}{{l+c-m-2-k}\choose l},
    \end{eqnarray*}                
where $\alpha'(i,l)$ defined in \eqref{dfn alph'}, $c_1=c-m-2$, $u_l= (-1)^l l!(2p+b-c-l)!(c-2)!$ and \begin{center}
    $v_i =\frac{(-1)^{i+c-m}(i-(p+b-c+2))!}{(i-1)! (m+1-i)! (c-m-2+i)!}.$
\end{center}
\end{lemma}

\begin{lemma}\label{app_C1_C2}
Let $n,c,i,l\in\mathbb{Z}$ such that $l\geq \text{max}\{0,n\}$ and $i\geq \text{max}\{1,n+1\}$. Also assume that $1\leq m\leq c-1$ and $n+c-m-1\geq 0$. Then
\begin{eqnarray*}
&\underset{0\leq k\leq c-m-1}\sum(-1)^k(c-k){{c-m-1+i}\choose k}{{c-m-2+i-k}\choose i-(n+1)}{{l+c-m-1-k}\choose l}\\
& = \begin{cases}
C_1+C_2 & \text{if}\ i\leq l\\
C_1+C_2+(-1)^{n+c-m-1+l}(m+1-i){{i-1}\choose l} & \text{if}\ l\leq i-1.
\end{cases}
\end{eqnarray*}
where 
$$C_1 =  \underset{0\leq j\leq  m_1}\sum (-1)^jc{{c-m-1+i}\choose i-(n+1+j)}{{l-(n+1+j)}\choose l-(n+c-m+j)}$$
$$C_2=\underset{0\leq j\leq m_2}\sum (-1)^j(n+c-m+j){{c-m-1+i}\choose i-(n+1+j)}{{l-(n+1+j)}\choose l-(n+c-m-1+j)}$$
where $m_1=\text{min}\{i-(n+1), l-(n+c-m)\}$ and $m_2=\text{min}\{i-(n+1), l-(n+c-m-1)\}$. Further, $C_1=0$ if $l<n+c-m$ and $C_2=0$ if $l<n+c-m-1$.
\end{lemma}
\begin{rk}
    When $m = c-1$ and $j =l-(n+c-m-1)$ then the binomial coefficient ${{l-(n+1+j)}\choose l-(n+c-m-1+j)} = {-1\choose 0} = 1$. This term appears in the last term of the sum in $C_2$ if $l-(n+c-m-1)\leq i-(n+1)$.
\end{rk}
\begin{proof}
First, we consider the following binomial expansion
$$(x-1)^{c-m-1+i}x^{-1} = \underset{0\leq k\leq c-m-1+i}\sum(-1)^k{{c-m-1+i}\choose k}x^{c-m-2+i-k}.$$
By differentiating $i-(n+1)$ times and dividing by $(i-n-1)!$, we get
\begin{align*}
&\hspace{1em}\underset{0\leq j\leq i-(n+1)}\sum(-1)^j{{c-m-1+i}\choose i-(n+1)-j}(x-1)^{n+c-m+j}x^{-(j+1)}\\
&=\underset{0\leq k\leq c-m-2+i}\sum(-1)^{k}{{c-m-1+i}\choose k}{{c-m-2+i-k}\choose i-(n+1)}x^{n+c-m-1-k}+(-1)^{n+c-m-2}x^{n-i}.
\end{align*}
Observe that ${{c-m-2+i-k}\choose i-(n+1)} = 0$ if $n+c-m-1<k\leq c-m-2+i$. Therefore, we have 
\begin{align*}
&\underset{0\leq j\leq i-(n+1)}\sum(-1)^j{{c-m-1+i}\choose i-(n+1)-j}(x-1)^{n+c-m+j}x^{-(j+1)}\\
&=\underset{0\leq k\leq n+c-m-1}\sum(-1)^{k}{{c-m-1+i}\choose k}{{c-m-2+i-k}\choose i-(n+1)}x^{n+c-m-1-k}+(-1)^{n+c-m-2}x^{n-i}.
\end{align*}
In the above equation, we multiply by $\frac{x^{l-n}}{l!}$, and then differentiate $l$ times. We get the following
\begin{align*}
&\underset{0\leq j\leq i-(n+1)}\sum(-1)^j{{c-m-1+i}\choose i-(n+1)-j}\left(\frac{D^l}{l!}\right)\left((x-1)^{n+c-m+j}x^{l-(n+1+j)}\right)
 \end{align*}
\begin{align*}
 &= \underset{0\leq k\leq n+c-m-1}\sum(-1)^k{{c-m-1+i}\choose k}{{c-m-2+i-k}\choose i-(n+1)}{{l+c-m-1-k}\choose l}x^{c-m-1-k}\\
& \hspace{22em}+\frac{(-1)^{n+c-m-2}}{l!}D^l(x^{l-i})
\end{align*}
where $D:=\frac{d}{dx}$ differential operator. Observe that ${{l+c-m-1-k}\choose l}=0$ if $c-m-1<k\leq n+c-m-1$ and ${{c-m-2+i-k}\choose i-(n+1)}=0$ if $n+c-m-1<k\leq c-m-1$. Therefore, even if $n <0$, we can extend the sum above to $k \leq c-m-1$, giving us 
\begin{align*}
&\underset{0\leq j\leq i-(n+1)}\sum(-1)^j{{c-m-1+i}\choose i-(n+1)-j}\left(\frac{D^l}{l!}\right)\left((x-1)^{n+c-m+j}x^{l-(n+1+j)}\right)\\
 &= \underset{0\leq k\leq c-m-1}\sum(-1)^k{{c-m-1+i}\choose k}{{c-m-2+i-k}\choose i-(n+1)}{{l+c-m-1-k}\choose l}x^{c-m-1-k}\\
& \hspace{22em}+\frac{(-1)^{n+c-m-2}}{l!}D^l(x^{l-i})
\end{align*}
 Again, we multiply by $x^{m+1}$ and after that differentiate one time the above equation to obtain the following 
\begin{align*}
&\underset{0\leq j\leq i-(n+1)}\sum(-1)^j{{c-m-1+i}\choose i-(n+1)-j}D\left(x^{m+1}\left(\frac{D^l}{l!}\right)\left((x-1)^{n+c-m+j}x^{l-(n+1+j)}\right)\right)\\
 &= \underset{0\leq k\leq c-m-1}\sum(-1)^k(c-k){{c-m-1+i}\choose k}{{c-m-2+i-k}\choose i-(n+1)}{{l+c-m-1-k}\choose l}x^{c-1-k}\\
& \hspace{22em}+\frac{(-1)^{n+c-m-2}}{l!}D(x^{m+1}D^l(x^{l-i}))
\end{align*}
By putting $x=1$ in the above equation, we get
\begin{align}\label{x=1_eqn}
\nonumber&\underset{0\leq k\leq c-m-1}\sum(-1)^k(c-k){{c-m-1+i}\choose k}{{c-m-2+i-k}\choose i-(n+1)}{{l+c-m-1-k}\choose l}\\
 \nonumber&=\underset{0\leq j\leq i-(n+1)}\sum(-1)^j{{c-m-1+i}\choose i-(n+1)-j}D\left(x^{m+1}\left(\frac{D^l}{l!}\right)\left((x-1)^{n+c-m+j}x^{l-(n+1+j)}\right)\right)|_{x=1}\\
 & \hspace{20em}+\frac{(-1)^{n+c-m-1}}{l!}D(x^{m+1}D^l(x^{l-i}))|_{x=1}.
 \end{align}
Now, 
\begin{align*}
&D\left(\frac{x^{m+1}D^l}{l!}\left((x-1)^{n+c-m+j}x^{l-(n+1+j)}\right)\right)\\
&=\underset{0\leq a\leq l}\sum{{n+c-m+j}\choose a}D\left((x-1)^{n+c-m+j-a}x^{m+1}\left(\frac{D^{l-a}}{(l-a)!}\right)(x^{l-(n+1+j)})\right).
\end{align*}
Observe that ${{n+c-m+j}\choose a} =0$ if $n+c-m+j< a$. Hence, we have
\begin{align*}
&D\left(\frac{x^{m+1}D^l}{l!}\left((x-1)^{n+c-m+j}x^{l-(n+1+j)}\right)\right)\\
&=\underset{0\leq a\leq n_1}\sum{{n+c-m+j}\choose a}D\left((x-1)^{n+c-m+j-a}x^{m+1}\left(\frac{D^{l-a}}{(l-a)!}\right)(x^{l-(n+1+j)})\right),
\end{align*}
where $n_1=\text{min}\{l, n+c-m+j\}$. Therefore, we get
\begin{align*}
&D\left(\frac{x^{m+1}D^l}{l!}\left((x-1)^{n+c-m+j}x^{l-(n+1+j)}\right)\right)\\
&=\underset{0\leq a\leq n_1}\sum{{n+c-m+j}\choose a}(n+c-m+j-a)\left((x-1)^{n+c-m-1+j-a}x^{m+1}\left(\frac{D^{l-a}}{(l-a)!}\right)(x^{l-(n+1+j)})\right)\\
&+\underset{0\leq a\leq n_1}\sum {{n+c-m+j}\choose a}\left((x-1)^{n+c-m+j-a}D\left(x^{m+1}\left(\frac{D^{l-a}}{(l-a)!}\right)(x^{l-(n+1+j)})\right)\right)
\end{align*}
Note that at $x=1$ the former sum is zero if $l<n+c-m-1$ and the latter sum is zero if $l<n+c-m$ (i.e., if $m_2<0$ and if $m_1<0$ respectively). We have
\begin{align*}
& \underset{0\leq j\leq i-(n+1)}\sum(-1)^j{{c-m-1+i}\choose i-(n+1)-j}D\left(x^{m+1}\left(\frac{D^l}{l!}\right)\left((x-1)^{n+c-m+j}x^{l-(n+1+j)}\right)\right)|_{x=1}\\
&= \underset{0\leq j\leq m_2}\sum (-1)^j(n+c-m+j){{c-m-1+i}\choose i-(n+1+j)}{{l-(n+1+j)}\choose l-(n+c-m-1+j)}\\
&\hspace{8em}+\underset{0\leq j\leq m_1}\sum(-1)^{j}c{{c-m-1+i}\choose i-(n+1+j)}{{l-(n+1+j)}\choose l-(n+c-m+j)}.
\end{align*}
Hence, from \eqref{x=1_eqn}, we get 
\begin{align*}
&\underset{0\leq k\leq c-m-1}\sum(-1)^k(c-k){{c-m-1+i}\choose k}{{c-m-2+i-k}\choose i-(n+1)}{{l+c-m-1-k}\choose l}\\
&=C_2+C_1+\frac{(-1)^{n+c-m-1}}{l!}D\left(x^{m+1}D^l(x^{l-i})\right)|_{x =1}.
\end{align*}
If $l\geq i$, then $D^l(x^{l-i}) = 0$ since  $l\geq i\geq 1$, giving us we the required identity in this case. If $l\leq i-1$, then
$$\frac{(-1)^{n+c-m-1}}{l!}D\left(x^{m+1}D^l(x^{l-i})\right)|_{x =1} = (-1)^{n+c-m-1+l}(m+1-i){{i-1}\choose l}.$$
Hence, we obtain our identity in this case too.
\end{proof}

\begin{lemma}\label{W_0_construction}
Assume all the hypotheses of Proposition \ref{rm_JH_conj}. Then there exists $W_0$ of Lemma \ref{lm m<c} for parts $(1)$ and $(3)$ in all cases, and for part $(2)$ only when $b\leq \nu$ of Proposition \ref{rm_JH_conj}.
\end{lemma}
\begin{proof}
We will prove our lemma in the following three cases.\\
\textbf{Case (1)}. $c-1\leq b\leq p$ and $\nu\leq c-1-\epsilon$\\
We begin by observing that we need to construct $W_0$  only when $b\geq c$ as $[0,n_1]$ is empty for $b =c-1$.\\
\textbf{Subcase (i)} $c\leq b\leq p-1$\\
We observe that $\nu\left({{r}\choose p-1}\right)\geq 1$ (as $c \geq 2$) and also ${r\choose b+j(p-1)}\equiv 0\bmod p$ for all $0\leq j\leq c-1$. Therefore, by Proposition \ref{gen 1} (for $l = 0$ and $m = 0$) we have $$(T-a_p)\left(\frac{f^0}{p}\right)\equiv\left[g,\ \underset{\substack{0<j<s\\j\equiv b\bmod (p-1)}}\sum \frac{{r\choose j}}{p}x^{r-j}y^j\right].$$
But $x^{r-j}y^j\equiv x^{r-\bar{j}}y^{\bar{j}}\bmod\left(V^{(1)}_r\right)$ where $\bar{j}\equiv j\bmod (p-1)$ and $2\leq \bar{j}\leq p$. Observe that $\frac{{r\choose j}}{p}\in\mathbb{Z}$ for all $0<j<s$ such that $j\equiv r\bmod (p-1)$. Hence, we have
$$ \underset{\substack{0<j<s \\ j\equiv r\bmod p-1}} \sum\frac{{r\choose j}}{p}x^{r-j}y^j\equiv \eta x^{r-b}y^b \mod \left( V^{(1)}_r\right)$$
where 
\begin{eqnarray*}
\eta & =  \underset{\substack{0<j<s\\ j\equiv s\bmod(p-1)}}\sum\frac{{r\choose j}}{p} 
\equiv  \underset{\substack{0<j<s\\ j\equiv s\bmod(p-1)}}\sum\frac{{s\choose j}}{p}
\equiv  \frac{b-s}{b}
\not\equiv  0\mod p.
\end{eqnarray*}
Here the first congruence follows as $\frac{{r\choose j}}{p}\equiv \frac{{s\choose j}}{p}\mod p^{t-1-\nu(j!)}$ and $t-1-\nu(j!)> \nu(a_p)>1$ (since $t>\nu(a_p)+c$, $\nu(j!)\leq \nu((s-(p-1))!)\leq c-1$). The second last congruence follows from Lemma $2.5$ of \cite{BG}. Hence, we have
$$(T-a_p)\left(\frac{f^0}{\eta p}\right)\equiv\left[g,\ x^{r-b}y^b\right]\bmod (V_r^{(1)})$$ 
as $\eta$ is $p$-adic unit. Therefore, there exists a $v_1\in V_r^{(1)}$ such that 
$$(T-a_p)\left(\frac{f^0}{\eta p}\right)\equiv\left[g,\ x^{r-b}y^b-v_1\right]$$ 
Using \eqref{r' not p} (with $r'=p-1+b$), we observe that the monomial $x^{r-b}y^b$ and $x^r$ generates the quotient $V_{p-1-b}\otimes D^b$ and submodule $V_b$ of $\frac{V_r}{V^{(1)}_r}$ respectively. But $x^r$ belongs to $\kerp$ by Remark 4.4 in \cite{KBG}.
Let $W_0$ is the sub module generated by $x^r$,\ $x^{r-b}y^b-v_1$ if $b\geq c$. Observe that $W_0$ satisfies all the required conditions.\\
\textbf{Subcase (ii)} $b = p$ \\
In this case by using (\ref{r' = p}) we have the following
$$0\longrightarrow V_1\longrightarrow\frac{V_r}{V^{(1)}_r}\longrightarrow V_{p-2}\otimes D\longrightarrow 0.$$
In the above exact sequence,  the first map sends $x$ to $x^r$ and the second map sends $x^{r-1}y$ to $x^{p-2}$. By the Remark $4.4$ of \cite{KBG}, we have $x^r,\ x^{r-1}y\in\kerp$ as $1 <\nu(a_p)$.  We define  $W_0$ in this case as the submodule generated by $x^r$ and $x^{r-1}y$, and observe that $W_0$ satisfies the required conditions of Lemma \ref{lm m<c}.

\textbf{Case (2)} $2\leq b\leq c-2$ and $1\leq \nu\leq \text{min}\{c-2,p+b-c\}$\\
In this case, we need to construct $W_0$ only when $b\leq \nu$ as it is clear from the statement of Proposition \ref{rm_JH_conj}.
By Proposition \ref{gen 1} (for $m = 0$ and $0\leq l\leq \text{min}\{\nu, p+b-c\}$), we have 
\begin{align*}
(T-a_p)\left(f^l\right)&\equiv\left[g,\ \underset{\substack{0\leq j\leq c-1}}\sum {{r-l}\choose b+ j(p-1)}x^{r-(b+j(p-1))}y^{(b+j(p-1))}\right]\\
&\equiv\left[g,\ \left(\underset{\substack{0\leq j\leq c-1}}\sum {{r-l}\choose b+ j(p-1)}\right)x^{r-b}y^b\right]\bmod (V_r^{(1)}).\\
\end{align*}
The last congruence follows from the same observation as in Case (1) above.  By applying Lemma $3.4$ in \cite{GK} for $m=0$ (see the second and fifth case), we get
\begin{align*}
{{r-l}\choose b+j(p-1)}
&\equiv \begin{cases}
{{p+b-c-l}\choose b-j}{{c-1}\choose j}\quad &\text{if}\ 0\leq l\leq p+b-c,\ 0\leq j\leq b\\
0\quad&\text{if}\ 0\leq l\leq p+b-c,\ b+1\leq j\leq c-1.\\
\end{cases}
\end{align*}
Therefore, we have 
\begin{align*}
(T-a_p)\left(f^l\right)&\equiv\left[g,\ \left(\underset{\substack{0\leq j\leq b}}\sum {{p+b-c-l}\choose b-j}{{c-1}\choose j}\right)x^{r-b}y^b\right]\bmod (V_r^{(1)})
\end{align*}
\begin{align*}
&\equiv\left[g,\ {{p+b-1-l}\choose b}x^{r-b}y^b\right]\bmod (V_r^{(1)}).
\end{align*}
The latter congruence is followed by Vandermonde's identity. Since $b\leq \nu$, one can take $l=b$ in Proposition \ref{gen 1}, in which case ${{p+b-1-l}\choose b}\not\equiv0(\bmod p)$. Thus, $$[g,x^{r-b}y^b]\in V_r^{(1)}+\kerp.$$ Using \eqref{r' not p} (with $r'=p-1+b$), we observe that the monomial $x^{r-b}y^b$ and $x^r$ generates the quotient $V_{p-1-b}\otimes D^b$ and submodule $V_b$ of $\frac{V_r}{V_r^{(1)}}$ respectively. But $x^r$ belongs to $\kerp$ by Remark 4.4 in \cite{KBG}. Therefore, taking $W_0$ to be the submodule generated by $x^r$ and $x^{r-b}x^b$ has the required properties.
                    
\textbf{Case (3)} $2\leq b\leq 2c-4-p$ and $p+b-c+1\leq\nu\leq c-3$\\
Observe that $b\leq \nu$ in this case, thus applying Proposition \ref{gen 1} with same ranges $m = 0$ and $0\leq l\leq \text{min}\{\nu, p+b-c\}$ as in Case $(2)$ giving required $W_0$.
\end{proof}


\bibliographystyle{plain}
\bibliography{LC_II_Ganguli_Kumar}
\end{document}